%% file: paper_SynRZME_02-06.tex
\documentclass[11pt,a4paper]{article}
\usepackage{authblk}
\usepackage{amsmath,amssymb}
\usepackage{graphicx}
\usepackage{array}
\usepackage{bm}
\usepackage{color}
\usepackage{multirow}
\usepackage{float}
\usepackage{enumerate}
\usepackage{caption}
\usepackage[super]{natbib}

\newcommand{\ud}{\mbox{d}}

\DeclareMathOperator*{\argmin}{arg\,min}

\usepackage{subfigure}
\usepackage{booktabs}
\usepackage[english]{babel}
\usepackage{multicol}
\usepackage[table]{xcolor}
\definecolor{light-gray}{gray}{0.9}
\newcolumntype{R}{!{\color{gray}\vline width 1pt}}

\title{Effect of the measurement errors on one-sided Synthetic-RZ control charts for monitoring the ratio of two normal variables}

\author[1,3]{Kim Duc Tran\thanks{Corresponding author. Email: \texttt{ductk@donga.edu.vn}}}
\author[2]{Thi Hien Nguyen}
\author[3]{Kim Phuc Tran}

\affil[1]{International Chair in DS \& XAI, International Research Institute for Artificial Intelligence and Data Science, Dong A University, Danang, Vietnam}

\affil[2]{Laboratoire AGM, UMR CNRS 8088, CY Cergy Paris Université, Cergy, France}

\affil[3]{Univ. Lille, ENSAIT, ULR 2461 -- GEMTEX -- Génie et Matériaux Textiles, F-59000 Lille, France}

\begin{document}

\maketitle

\begin{abstract}
In numerous industrial production settings, keeping track of the ratio formed by two normally distributed random variables is a task of considerable practical interest. The present work examines how measurement errors influence the behaviour of a pair of one-sided Synthetic control charts designed to monitor such a ratio (referred to here as Synthetic-RZ charts), with the analysis covering both the zero-state and the steady-state average run length ($ARL$). To incorporate measurement error into the operation of these charts, we adopt a linear covariate error model. We describe, step by step, how the parameters of the underlying model evolve as the process moves from an in-control to an out-of-control state, and we deliberately avoid the restrictive premise that the observed shift magnitude is unrelated to the measurement errors. The run length characteristics of the charts are obtained by means of a Markov chain formulation. A series of numerical experiments makes clear that measurement error erodes the detection capability of the charts. A particularly useful outcome of the investigation is that collecting several measurements on each inspected unit does not constitute an efficient remedy for the adverse influence of measurement error on the performance of the Synthetic-RZ charts.
\end{abstract}


\section{Introduction}
\label{sec:introduction}
Statistical Process Control (SPC) has long served as an effective methodology for overseeing manufacturing operations in industry. Within the SPC toolkit, control charts occupy a central position because of their ability to flag departures of a process from its intended behaviour. A large body of work has been devoted to constructing charts that respond more rapidly to process shifts; representative contributions include \citet{Tran2017Runrules_median_mean}, \citet{TranPH-SyntheticME}, \citet{Yeong_CV_ME_2017}, and \citet{Nguyen-Ratio-Me}. The earliest and arguably most popular scheme is the Shewhart chart, proposed by \citet{Shewhart1931}, whose appeal lies chiefly in how straightforward it is to apply. That said, when shifts are small or of moderate magnitude, the Shewhart chart reacts more sluggishly than memory-based alternatives such as the CUSUM chart \citet{Brook1972} and the EWMA chart \citet{Hunter1986}. The Synthetic chart, introduced by \citet{Wu2000}, was conceived to retain the simplicity of the Shewhart philosophy while improving sensitivity. It operates by coupling a Shewhart sub-chart with a conforming run length ($CRL$) sub-chart, and has been reported to compete favourably with rival schemes when detecting a range of process shifts.

In a variety of industrial production contexts, the quantity that genuinely matters is the ratio of two normally distributed variables, and several concrete instances of this have been documented. \citet{Tran2016_Shewhart_RZ_ME} described a case drawn from battery-recycling facilities in Italy, in which the ratio of the weight of ``recyclable batteries'' to that of the ``total batch'' must be tracked so as to gauge the associated economic loss. A second illustration, from the food sector, was given by \citet{Celano2016_phaseII_RZ}, where the quantity of interest is the weight ratio of ``pumpkin seeds'' to ``flaxseeds''. Further manufacturing situations of this kind are noted in \citet{Tran2016_Runrules_RZ,Tran2016_CUSUM_RZ}. A range of charts has accordingly been put forward for monitoring such a ratio, among them \citet{Nguyen_VSIEWMARZ_2018,Tran2017_steadystate,Tran2016_Runrules_RZ,Tran2015_EWMA_RZ,Tran2016_CUSUM_RZ}. A Synthetic chart for this purpose---the Synthetic-RZ chart---was developed by \citet{Celano2016_Synthentic_RZ}. Their proposal, however, is two-sided, and because the distribution of the ratio of two normal variables is skewed, the chart turns out to be $ARL$-biased: situations arise in which the out-of-control $ARL$ exceeds the in-control one. As a remedy, \citet{NICS19} advocated running two distinct one-sided Synthetic charts, one tuned to detect a downward move in the ratio and the other an upward move. Relative to the two-sided design of \citet{Celano2016_Synthentic_RZ}, this pair of one-sided charts removes the $ARL$-bias and detects ratio shifts more effectively.

With a view to making ratio charts more usable in practice, attention has recently turned to the question of measurement error. The topic itself has a substantial literature---see, for instance, \citet{Linna2001b}, \citet{Tran2017ShewhartMedianME}, \citet{Tran_2018_CVME}, \citet{Nguyen_Thong_Chatti_IEEE2018}, \citet{Tranhanh2016_CUSUM_CV}, \citet{Yeong_CV_ME_2017}, \citet{Cheng2017_EWMA_median_with_ME}, and \citet{Nguyen_2019_VSICV_ME}---and these works collectively demonstrate that measurement error exerts a pronounced influence both on the control limits and on chart performance. In the specific setting of ratio charts, \citet{Tran2016_Shewhart_RZ_ME} were the first to examine the impact of measurement error, doing so for a Shewhart-RZ chart. Their analysis, however, rested on a fairly restrictive premise, namely that the observed process shift is independent of the measurement error. One contribution of the present paper is to dispense with that premise so as to obtain a more realistic treatment. Furthermore, while the effect of measurement error has by now been studied for the Shewhart-RZ and the EWMA-RZ charts, the Synthetic-RZ charts---which tend to strike a useful compromise between simplicity and detection ability---have not yet been investigated in this respect. We therefore study how measurement error affects the two one-sided Synthetic-RZ charts of \citet{NICS19}, representing the error through a linear covariate error model and deriving in detail how the process parameters change in the presence of error, again without invoking independence between the shift magnitude and the measurement errors.

The remainder of the paper proceeds as follows. Section \ref{sec:distribution} recalls briefly the sampling distribution of the ratio of two normal variables. Section \ref{sec:errormodel} sets out a linear covariate error model for the ratio and establishes how the model parameters move from the in-control to the out-of-control state. Section \ref{sec:implementation} addresses the design and operation of the one-sided Synthetic-RZ charts when measurement error is present, and Section \ref{sec:numerical} reports the effect of measurement error on their performance. A worked example illustrating the use of one such chart under measurement error is given in Section \ref{sec:illustrative}, and Section \ref{sec:conclusions} offers some closing comments.

\section{A short overview of the sample distribution of the ratio}
\label{sec:distribution}
Let $\mathbf{W}=(X,Y)^\intercal$ be a bivariate normal random vector with mean vector $\boldsymbol{\mu_W}$ and variance--covariance matrix $\boldsymbol{\Sigma_W}$ given by
\begin{equation}
\label{equ:variance}
\boldsymbol{\mu}_W=
\begin{bmatrix}
\mu_X \\
\mu_Y
\end{bmatrix}
\qquad \text{and} \qquad
\boldsymbol{\Sigma}_W=
\begin{bmatrix}
\sigma_X^2 & \rho \sigma_X \sigma_Y \\
\rho \sigma_X \sigma_Y & \sigma_Y^2
\end{bmatrix}.
\end{equation}
Here $\rho$ denotes the correlation between $X$ and $Y$, while $\sigma_X$ and $\sigma_Y$ are their respective standard deviations and $\mu_X$, $\mu_Y$ their means. The two coefficients of variation and the ratio of the standard deviations are then $\gamma_X=\frac{\sigma_X}{\mu_X}$, $\gamma_Y=\frac{\sigma_Y}{\mu_Y}$ and $\omega=\frac{\sigma_X}{\sigma_Y}$, respectively.

The ratio of $X$ to $Y$ is written $Z=X/Y$. The distribution of $Z$ has been the subject of several investigations---e.g. \citet{Geary1930}, \citet{Hayya1975}, \citet{Pham2006}. \citet{Cedilnik2004} derived a general closed-form expression for its probability density function ($p.d.f$), and its cumulative distribution function ($c.d.f$) along with the inverse $c.d.f$ can be evaluated numerically following \citet{Celano2014a}. Since these exact routes are somewhat cumbersome to apply, approximations are usually preferred. As shown by \citet{Celano2016_phaseII_RZ}, provided the coefficients of variation of $X$ and $Y$ are not too large, $Z$ admits the accurate approximation
\begin{equation}
\label{equ:CDFZ}
F_Z(z|\gamma_X,\gamma_Y,\omega,\rho)\simeq\Phi\left(\frac{A}{B}\right),
\end{equation}
\noindent
with $\Phi(\cdot)$ the $c.d.f$ of the standard normal distribution, and where the quantities $A$ and $B$ depend on $z$, $\gamma_X$, $\gamma_Y$, $\omega$ and $\rho$ through
\begin{eqnarray*}
A & = & \frac{z}{\gamma_Y}-\frac{\omega}{\gamma_X}, \\
B & = & \sqrt{\omega^2-2\rho\omega z+z^2}.
\end{eqnarray*}
For a stable process whose quality variables are normally distributed, the spread of the population ought to be much smaller than its mean, so assuming small $\gamma_X$ and $\gamma_Y$---and hence adopting the approximate distribution---is justified. Differentiating the $c.d.f$ in \eqref{equ:CDFZ} gives the $p.d.f$ of $Z$ as
\begin{equation}
\label{equ:PDFZ}
f_Z(z|\gamma_X,\gamma_Y,\omega,\rho)\simeq\left(\frac{1}{B\gamma_Y}
-\frac{(z-\rho\omega)A}{B^3}\right)\times\phi\left(\frac{A}{B}\right),
\end{equation}
with $\phi(\cdot)$ the standard normal density.

An approximation to the inverse distribution function ($i.d.f$) $F^{-1}_Z(p|\gamma_X,\gamma_Y,\omega,\rho)$ of $Z$ can be obtained in the same spirit, namely
\begin{equation}
\label{equ:IDFZ}
F^{-1}_Z(p|\gamma_X,\gamma_Y,\omega,\rho)\simeq\left\{
\begin{array}{ll}
\frac{-C_2-\sqrt{C_2^2-4C_1C_3}}{2C_1} & \mbox{if }p\in(0,0.5], \\
\frac{-C_2+\sqrt{C_2^2-4C_1C_3}}{2C_1} & \mbox{if }p\in[0.5,1),
\end{array}
\right.
\end{equation}
\noindent
where the coefficients $C_1$, $C_2$ and $C_3$ depend on $p$, $\gamma_X$, $\gamma_Y$, $\omega$ and $\rho$ via
\begin{eqnarray*}
C_1 & = & \frac{1}{\gamma_Y^2}-\Phi^{-1}(p)^2, \\
C_2 & = & 2\omega\left(\rho \Phi^{-1}(p)^2-\frac{1}{\gamma_X\gamma_Y}\right), \\
C_3 & = & \omega^2\left(\frac{1}{\gamma_X^2}-\Phi^{-1}(p)^2\right),
\end{eqnarray*}
\noindent
and $\Phi^{-1}(\cdot)$ is the $i.d.f$ of the standard normal distribution.

\section{Linear covariate error model for the sample of the ratio}
\label{sec:errormodel}
This section develops a linear covariate error model for the sample ratio and gives a precise account of how the process parameters are altered once a shift occurs in the presence of measurement error.

Return to the bivariate normal vector $\mathbf{W}\sim N(\boldsymbol{\mu}_{\mathbf{W}},\boldsymbol{\Sigma}_{\mathbf{W}})$ introduced in \eqref{equ:variance}. Suppose $\{\mathbf{W}_{i,1},\mathbf{W}_{i,2},\ldots,\mathbf{W}_{i,n}\}$ is a collection of $n$ independent draws from it, where $\mathbf{W}_{i,j}=(X_{i,j},Y_{i,j})^\intercal$ is the quality characteristic of the $j$-th unit, $j=1,\ldots,n$, inspected at time $i$, $i=1,2,\ldots$

Owing to measurement error, the \emph{true} characteristic $\mathbf{W}_{i,j}$ cannot be observed directly; what is available instead is a set of $m\geqslant 1$ repeated readings on unit $j$ at time $i$, say $\{\mathbf{W}^*_{i,j,1},\mathbf{W}^*_{i,j,2},\ldots,\mathbf{W}^*_{i,j,m}\}$. Following the linear covariate error model of \citet{LinnaWoodall2001}, each reading relates to the true value through
\begin{equation}
\label{equ:model_ME}
\mathbf{W}^*_{i,j,k}=\mathbf{A}+\mathbf{B}\mathbf{W}_{i,j}+\boldsymbol{\varepsilon}_{i,j,k}, \quad k=1,\ldots,m,
\end{equation}
\noindent
in which $\mathbf{A} = (a_X,a_Y)^\intercal$ is a $(2\times 1)$ vector of constants, $\mathbf{B}$ is a $(2\times 2)$ matrix and $\boldsymbol{\varepsilon}_{i,j,k}\sim N(\mathbf{0},\boldsymbol{\Sigma}_M)$ is a bivariate normal vector independent of $\mathbf{W}_{i,j}$. Denoting by $\sigma_{MX}$ and $\sigma_{MY}$ the standard deviations of the measurement error in $X$ and $Y$ and by $\rho_M\in(-1,+1)$ their correlation, the covariance matrix of $\boldsymbol{\varepsilon}$ reads
\begin{equation}
\boldsymbol{\Sigma}_M=\left(
\begin{array}{cc}
\sigma_{MX}^2 & \rho_M\sigma_{MX}\sigma_{MY} \\
\rho_M\sigma_{MX}\sigma_{MY} & \sigma_{MY}^2
\end{array}\right).
\end{equation}
As in \citet{Tran2016_Shewhart_RZ_ME}, we take $\mathbf{B}=\mathbf{I}_{2\times2}$, the identity matrix, which amounts to assuming a constant accuracy-error vector $\mathbf{A}$ across the measurement range.

Because $\mathbf{W}_{i,j}$ itself is hidden, its sample mean over the repeated readings,
\begin{eqnarray}
\label{equ:ME model}
\overline{\mathbf{W}}^*_{i,j} = \frac{1}{m}\sum_{k=1}^m \mathbf{W}^*_{i,j,k} = \mathbf{A}+\mathbf{B}\mathbf{W}_{i,j}+\frac{1}{m}\sum_{k=1}^m\boldsymbol{\varepsilon}_{i,j,k},
\end{eqnarray}
is commonly used in its place. Writing $\overline{\mathbf{W}}^*_{i,j}=(\bar{X}^*_{i,j},\bar{Y}^*_{i,j})$, it follows from \eqref{equ:ME model} that this average is again a bivariate normal vector, with mean vector
\begin{eqnarray}
\label{equ:mustar}
\boldsymbol{\mu}_{\mathbf{W}^*} & = & \mathbf{A}+\mathbf{B}\boldsymbol{\mu_W}
\end{eqnarray}
and variance--covariance matrix
\begin{eqnarray}
\boldsymbol{\Sigma}_{\mathbf{W}^*} & = & \mathbf{B}\boldsymbol{\Sigma_W}\mathbf{B}^\intercal + \frac{1}{m}\boldsymbol{\Sigma}_M =
\boldsymbol{\Sigma_W} + \frac{1}{m}\boldsymbol{\Sigma}_M.
\label{equ:sigmastar}
\end{eqnarray}

Building on these expressions, we now trace how the process parameters respond to a shift under measurement error.

While the process is in control, write the mean vector of the true characteristic $\mathbf{W}_{i,j} = (X_{i,j},Y_{i,j})$ as $\boldsymbol{\mu_{0,W}}=\left( \mu_{0,X},\mu_{0,Y} \right)^\intercal$ and the correlation between $X_{i,j}$ and $Y_{i,j}$ as $\rho=\rho_0$, so that the mean ratio equals $z_0=\frac{\mu_{0,X}}{\mu_{0,Y}}$. Suppose that, when the process becomes disturbed, the ratio moves from $z_0$ to $z_1=\tau z_0$ (with $\tau$ the shift size) and the correlation moves from $\rho_0$ to $\rho_1$. The change in $z_0$ corresponds to a displacement of $\boldsymbol{\mu_{W}}$ from $\boldsymbol{\mu_{0,W}}$ to $\boldsymbol{\mu_{1,W}}=\left( \mu_{0,X}+\delta_X\sigma_X, \mu_{0,Y}+\delta_Y\sigma_Y\right)^\intercal$, where $\delta_X$ and $\delta_Y$ measure the size of the mean shifts in $X_{i,j}$ and $Y_{i,j}$ and the two standard deviations $\sigma_X$ and $\sigma_Y$ are assumed equal. The out-of-control mean ratio is therefore
\[
z_1= \frac{\mu_{0,X}+\delta_X \sigma_X }{\mu_{0,Y}+\delta_Y \sigma_Y }= \tau \times z_0=\tau\times \frac{\mu_{0,X}}{\mu_{0,Y}},
\]
from which the shift size satisfies
\[\tau=\frac{1+\delta_X \gamma_X}{1+\delta_Y \gamma_Y},\]
or equivalently
\begin{equation}
\label{equ:tau}
1+\delta_X \gamma_X =\tau(1+\delta_Y \gamma_Y).
\end{equation}

Writing the mean vector and covariance matrix of $\overline{\mathbf{W}}^*_{i,j}$ as
\begin{equation}
\label{equ:variance ME}
\boldsymbol{\mu_{W^*}}=\left(\begin{array}{c} \mu_{X^*} \\ \mu_{Y^*} \end{array}\right)~ \text{and} ~ \boldsymbol{\Sigma_{W^*}}=\left(\begin{array}{cc}
\sigma_{X^*}^2 & \rho^*\sigma_{X^*}\sigma_{Y^*} \\
\rho^*\sigma_{X^*}\sigma_{Y^*} & \sigma_{Y^*}^2 \end{array} \right),
\end{equation}
equations \eqref{equ:mustar} and \eqref{equ:sigmastar} yield
\begin{eqnarray}
\label{mustarx}
\mu_{X^*}&= &a_X+\mu_X+\delta_X\sigma_X,\\
\mu_{Y^*}&=& a_Y+\mu_Y+\delta_Y\sigma_Y,\\
\label{sigmastarx}
\sigma_{X^*}^2 & = & \sigma^2_X+\frac{\sigma_{MX}^2}{m}, \\
\label{sigmastary}
\sigma_{Y^*}^2 & = & \sigma^2_Y+\frac{\sigma_{MY}^2}{m},\\
\label{rhostar}
\rho^*&=&\frac{\rho \sigma_X\sigma_Y +\rho_M\frac{\sigma_{MX}\sigma_{MY}}{m}}
{\sigma_{X^*}\sigma_{Y^*}}.
\end{eqnarray}
Hence the coefficients of variation $\gamma_{X^*}=\frac{\sigma_{X^*}}{\mu_{X^*}}$ and $\gamma_{Y^*}=\frac{\sigma_{Y^*}}{\mu_{Y^*}}$ of $\bar{X}_{i,j}^*$ and $\bar{Y}_{i,j}^*$ become
\begin{eqnarray}
\label{equ:gammaxo}
\gamma_{X^*} & = & \frac{\sqrt{\sigma^2_X+\frac{\sigma_{MX}^2}{m}}}
{a_X+ \mu_{0,X}+\delta_X\sigma_X }, \\
\label{equ:gammayo}
\gamma_{Y^*} & = & \frac{\sqrt{\sigma^2_Y+\frac{\sigma_{MY}^2}{m}}}
{a_Y+ \mu_{0,Y}+\delta_Y\sigma_Y}.
\end{eqnarray}
Dividing the numerator of \eqref{equ:gammaxo} by $\sigma_X$ and the denominator by $\mu_{0,X}$, and then substituting $\tau(1+\delta_Y\gamma_Y)$ for $1+\delta_X\gamma_X$ according to \eqref{equ:tau}, recasts $\gamma_{X^*}$ as
\begin{eqnarray}
\label{equ:gammaxf}
\gamma_{X^*}&=&\frac{\sqrt{1+\frac{\eta^2_X}{m}}}{1+\delta_X \gamma_X +\theta_X}
\times\gamma_X=\frac{\sqrt{1+\frac{\eta^2_X}{m}}}{\tau(1+\delta_Y \gamma_Y) +\theta_X}
\times\gamma_X,
\end{eqnarray}
where $\eta_X=\frac{\sigma_{MX}}{\sigma_X}$, $\theta_X=\frac{a_X}{\mu_{0,X}}$ and $\gamma_X=\frac{\sigma_X}{\mu_{0,X}}$.

In the same manner, $\gamma_{Y^*}$ from \eqref{equ:gammayo} and $\rho^*$ from \eqref{rhostar} take the forms
\begin{eqnarray}
\label{equ:gammayf}
\gamma_{Y^*}&=& \frac{\sqrt{1+\frac{\eta^2_Y}{m}}}{1+\delta_Y \gamma_Y +\theta_Y}\times\gamma_Y,\\
\label{equ:rhof}
\rho^*&=&\frac{\rho+\rho_M\frac{\eta_X\eta_Y}{m}}
{\sqrt{1+\eta^2_X/m}\sqrt{1+\eta^2_Y/m}},
\end{eqnarray}
with $\eta_Y=\frac{\sigma_{MY}}{\sigma_Y}$, $\theta_Y=\frac{a_Y}{\mu_{0,Y}}$ and $\gamma_Y=\frac{\sigma_Y}{\mu_{0,Y}}$.

Likewise, the ratio of standard deviations $\omega^*=\frac{\sigma_{X^*}}{\sigma_{Y^*}}$ equals
\begin{eqnarray}
\label{equ:omegaf}
\omega^*&=&\sqrt{\frac{1+\frac{\eta^2_X}{m}}{1 +\frac{\eta^2_Y}{m}}} \times\omega,
\end{eqnarray}
where $\omega=\frac{\sigma_X}{\sigma_Y}$.

With this notation, the in-control and out-of-control ratios in the presence of measurement error are
\begin{eqnarray}
\label{equ:z0s}
z^*_0&=&\frac{\mu_{0,X^*}}{\mu_{0,Y^*}}= \frac{\mu_{0,X}+a_X }{\mu_{0,Y}+a_Y }= \frac{1+\theta_X}
{1+\theta_Y}\times z_0,\\
\label{equ:z1s}
z^*_1&=&\frac{\mu_{1,X^*}}{\mu_{1,Y^*}}= \frac{a_X+\mu_{0,X}+\delta_X \sigma_X }{a_Y+ \mu_{0,Y}+\delta_Y \sigma_Y }= \frac{1+\theta_X+\delta_X \gamma_X}
{1+\theta_Y+\delta_Y \gamma_Y}\times z_0.
\end{eqnarray}
Equations \eqref{equ:z0s}--\eqref{equ:z1s} make clear that, in general, $z_1^*\neq \tau z_0^*$. This is precisely where our treatment departs from \citet{Tran2016_Shewhart_RZ_ME}, who assumed the observed shift size to be independent of the measurement error; abandoning that assumption affords a more realistic description of how measurement error bears on ratio control charts.

\section{Design and implementation of the Synthetic-RZ control chart with measurement error}
\label{sec:implementation}
\subsection{Distribution of the monitored statistic}

When measurement error is present, the statistic being monitored takes the form
\begin{equation}\label{equa:RZ}
\hat{Z}^*_i=\frac{\hat{\mu}_{X_i^*}}{\hat{\mu}_{Y_i^*}}=\frac{\bar{X}^*_i}{\bar{Y}^*_i},
\end{equation}
with $\bar{X}^*_i=\frac{1}{n}\sum_{j=1}^n \bar{X}^*_{i,j}$ and $\bar{Y}^*_i=\frac{1}{n}\sum_{j=1}^n \bar{Y}^*_{i,j}$, where $\bar{X}^*_{i,j}$ and $\bar{Y}^*_{i,j}$ are the components of the vector $\overline{\mathbf{W}}^*_{i,j}$ of \eqref{equ:ME model}.

By construction, $\bar{X}^*_i\sim N(\mu_{X^*},\frac{\sigma_{X^*}}{\sqrt{n}})$ and $\bar{Y}^*_i\sim N(\mu_{Y^*},\frac{\sigma_{Y^*}}{\sqrt{n}})$, so the coefficients of variation of $\bar{X}^*_i$ and $\bar{Y}^*_i$ are
\begin{eqnarray}
\label{equ:gammaxsm}
\gamma_{\bar{X}^*} & = & \frac{\sigma_{X^*}}{\mu_{X^*}\sqrt{n}}
=\frac{\gamma_{X^*}}{\sqrt{n}}, \\
\label{equ:gammaysm}
\gamma_{\bar{Y}^*} & = & \frac{\sigma_{Y^*}}{\mu_{Y^*}\sqrt{n}}
=\frac{\gamma_{Y^*}}{\sqrt{n}},
\end{eqnarray}
while the ratio of their standard deviations $\omega_i^*$ is
\begin{equation}
\label{equ:omegasm}
\omega_i^*=\frac{\sigma_{X^*}/\sqrt{n}}{\sigma_{Y^*}/\sqrt{n}}
=\frac{\sigma_{X^*}}{\sigma_{Y^*}}=\omega^*.
\end{equation}

The $c.d.f$ and the $i.d.f$ of $\hat{Z}^*_i$ can accordingly be written as
\begin{eqnarray}
\label{cdfZi}
F_{\hat{Z}^*_i}(z\mid n,\gamma_{X^*},\gamma_{Y^*},z^*_0,\rho^*_0)&=&F_{Z^*}\left(z\mid \frac{\gamma_{X^*}}{\sqrt{n}},\frac{\gamma_{Y^*}}{\sqrt{n}},\frac{z^*_0\gamma_{X^*}}{\gamma_{Y^*}},\rho^*_0\right),\\
\label{idfZi}
F_{\hat{Z}^*_i}^{-1}(p\mid n,\gamma_{X^*},\gamma_{Y^*},z^*_0,\rho^*_0)&=&F_{Z^*}^{-1}\left(p\mid \frac{\gamma_{X^*}}{\sqrt{n}},\frac{\gamma_{Y^*}}{\sqrt{n}},\frac{z^*_0\gamma_{X^*}}{\gamma_{Y^*}},\rho^*_0\right),
\end{eqnarray}
where $F_{Z^*}(\ldots)$ and $F_{Z^*}^{-1}(\ldots)$ are those of \eqref{equ:CDFZ} and \eqref{equ:IDFZ}, evaluated with the parameters $\gamma_{X^*}$, $\gamma_{Y^*}$, $\omega^*$ and $\rho^*$ supplied by \eqref{equ:gammaxf}, \eqref{equ:gammayf}, \eqref{equ:omegaf} and \eqref{equ:rhof}.

\subsection{Design of the one-sided Synthetic-RZ control charts}

A Synthetic chart pairs a Shewhart sub-chart with a conforming run length ($CRL$) sub-chart. The $CRL$ counts the samples gathered between one nonconforming sample and the next; when no preceding nonconforming sample exists, it counts from the start of monitoring up to the first nonconforming sample \citet{Costa2006}.

As noted in Section \ref{sec:introduction}, the asymmetry of the distribution of $Z$ (cf. \eqref{equ:CDFZ}--\eqref{equ:PDFZ}) means that a single two-sided Synthetic chart watching for both an increase and a decrease in $Z$ is $ARL$-biased, as in \citet{Celano2016_Synthentic_RZ}. Following \citet{NICS19}, we sidestep this by employing two separate one-sided charts. The first is an upper-sided chart aimed at spotting an increase in the ratio (hereafter the `Synthetic-RZ$^+$ chart'); the second is a lower-sided chart aimed at spotting a decrease (the `Synthetic-RZ$^-$ chart').

Designing the one-sided Synthetic-RZ-ME charts means fixing, for the lower-sided chart, the lower control limit $LCL^-$ of its Shewhart-RZ sub-chart together with the limit $H^-$ of its $CRL$ sub-chart, and, for the upper-sided chart, the upper control limit $UCL^+$ together with the limit $H^+$. The two charts are run according to the steps below:

\begin{itemize}
\item \textbf{Step 1.} Fix the target pair $(H^-, LCL^-)$ for the lower-sided chart (or $(H^+, UCL^+)$ for the upper-sided one).

\item \textbf{Step 2.} At time $i$ ($i=1,2,\ldots$), draw a pair of random samples of size $n\times m$, written $(X^*_{i,j,k},Y^*_{i,j,k})$, and form the statistic $\hat{Z}^*_i$ via \eqref{equa:RZ}.

\item \textbf{Step 3.} If the sample conforms---$\hat{Z}^*_i < UCL^+$ for the upper-sided chart, or $\hat{Z}^*_i>LCL^-$ for the lower-sided chart---go back to Step 2 with time advanced to $i+1$; otherwise continue to Step 4.

\item \textbf{Step 4.} Evaluate the $CRL$ and compare it against the $CRL$ sub-chart limit ($H^+$ or $H^-$). If $CRL>H^+$ (or $CRL>H^-$), treat the process as in control and return to Step 2; otherwise proceed to Step 5.

\item \textbf{Step 5.} Declare the process out of control, locate and eliminate the assignable cause, and then resume at Step 2.
\end{itemize}

\subsection{Run length properties via the Markov chain approach}

Chart performance is customarily judged by the average run length ($ARL$), the expected number of samples collected before the first out-of-control alarm; a smaller $ARL$ signifies better detection. Our aim is thus to choose the control limits $(H^{-*},LCL^{-*})$ or $(H^{+*},UCL^{+*})$ so that, for a prescribed out-of-control scenario, the $ARL$ is minimised while the in-control $ARL$ is held at a preset level $ARL_0$. The two one-sided design problems can be stated as follows:
\begin{itemize}
\item for the Synthetic-RZ$^-$-ME chart,
\begin{align}\label{eq:low}
(H^{-*},LCL^{-*})=\argmin_{(H^{-},LCL^-)} ARL(H^{-}, LCL^-,n,\gamma_{X^*},\gamma_{Y^*},\rho^*_0, \rho^*_1,\tau)
\end{align}
subject to
\begin{equation}
\label{equ:ARL0low}
ARL(H^{-}, LCL^-,n,\gamma_{X^*},\gamma_{Y^*},\rho^*_0, \tau=1)=ARL_0,
\end{equation}
\item for the Synthetic-RZ$^+$-ME chart,
\begin{align}\label{eq:up}
(H^{+*},UCL^{+*})=\argmin_{(H^{+},UCL^+)} ARL(H^{+}, UCL^+,n,\gamma_{X^*},\gamma_{Y^*},\rho^*_0, \rho^*_1,\tau)
\end{align}
subject to
\begin{equation}
\label{equ:ARL0up}
ARL(H^{+}, UCL^+,n,\gamma_{X^*},\gamma_{Y^*},\rho^*_0, \tau=1)=ARL_0.
\end{equation}
\end{itemize}

We describe below how the $ARL$ of the Synthetic-RZ$^-$-ME chart is obtained; the Synthetic-RZ$^+$-ME case follows along the same lines.

To derive the run-length properties of the Synthetic-RZ$^-$-ME chart we adopt a Markov chain whose transition probability matrix $\mathbf{P}$ is
\begin{eqnarray*}
\mathbf{P}=
\begin{pmatrix}
\mathbf{Q} & \mathbf{r} \\
\mathbf{0}^{\intercal} & 1
\end{pmatrix}
=
\left(
\begin{array}{*{7}c}
1-\theta &\theta & 0 & \ldots & 0 & \vline & 0 \\
0 & 0 & 1-\theta & & 0 & \vline & \theta \\
\vdots & & & \ddots & \vdots & \vline & \vdots \\
0 &\ldots& \ldots & 0 &1-\theta & \vline & \theta \\
1-\theta& 0 & \ldots & \ldots & 0 & \vline & \theta \\
\hline
0 & \ldots & \ldots & \ldots & 0 & \vline & 1 \\
\end{array}
\right),
\end{eqnarray*}
in which $\mathbf{Q}$ is the $(H^{-}+1)\times(H^{-}+1)$ sub-matrix governing transitions among the transient states, $\mathbf{0}^{\intercal}=(0,0,\dots,0)$ is a $1\times(H^{-}+1)$ row vector, $\mathbf{r}$ is the $(H^{-}+1)\times 1$ column vector with $\mathbf{r}=\mathbf{1}-\mathbf{Q}\mathbf{1}$ and $\mathbf{1}=(1,1,\dots,1)^{\intercal}$, and $\theta$ is the probability that a sample is nonconforming. Specifically,
\begin{itemize}
\item for the downward chart,
\begin{eqnarray}
\label{equ:theta down}
\theta&=& F_{\hat{Z}^*_i}(LCL^-\mid n,\gamma_{X^*},\gamma_{Y^*},z^*_1,\rho^*_1);
\end{eqnarray}
\item for the upward chart,
\begin{eqnarray}
\label{equ:theta up}
\theta&=& 1-F_{\hat{Z}^*_i}(UCL^+\mid n,\gamma_{X^*},\gamma_{Y^*},z^*_1,\rho^*_1),
\end{eqnarray}
\end{itemize}
with $F_{\hat{Z}^*_i}(\ldots)$ as in \eqref{cdfZi}.

Letting $\mathbf{q}=(q_0,q_1,\ldots,q_{H^{-}})^{\intercal}$ denote the $(H^{-}+1)\times 1$ vector of initial probabilities over the transient states, the zero-state $ARL$ is given by (see \citet{Davis2002})
\begin{eqnarray}
\label{equ:ARL-zero}
ARL = \mathbf{q}^{\intercal}\mathbf{I}-\mathbf{Q}^{-1}\mathbf{1}.
\end{eqnarray}

In practice a process that has run in control for some time settles into a steady state, and the long-run behaviour of the chart is then better captured by the steady-state $ARL$. The Markov chain again furnishes the \textit{cyclical} steady-state $ARL$:
\begin{equation}
\label{equ:ARL-steady state}
ARL = \boldsymbol{\psi}^{\intercal} \mathbf{I}-\mathbf{Q}^{-1}\mathbf{1},
\end{equation}
in which the \textit{cyclical} steady-state vector $\boldsymbol{\psi}$ is
$\boldsymbol{\psi}
=
\dfrac{\mathbf{I}-(\mathbf{Q}^{\intercal})^{-1}\mathbf{q}}
{\mathbf{1}^{\intercal}\mathbf{I}
-
\mathbf{1}^{\intercal}(\mathbf{Q}^{\intercal})^{-1}\mathbf{q}}$,
as established by \citet{Darroch1965}.

Frequently the exact shift size $\tau$ is unknown when the chart is designed, and in that event relying on the $ARL$ as a performance measure may produce a poorly tuned chart. For this reason we also gauge performance by the expected average run length ($EARL$),
\begin{equation}
\label{equ:EARL}
EARL=\int_{\Omega}ARL\times f_{\tau}(\tau)\ud\tau,
\end{equation}
where $f_{\tau}(\tau)$ is the density of the random shift $\tau$ over a support $\Omega$. In keeping with common SPC practice, $\tau$ is taken to be uniformly distributed on a chosen interval $[a,b]$ (see, for example, \citet{Nguyen_VSIEWMARZ_2018}), so that $f_{\tau}(\tau)=\frac{1}{b-a}$ for $\tau\in\Omega=[a,b]$. As in \citet{Tran2016_Shewhart_RZ_ME}, two ranges are considered: $\Omega_D=[0.9,1)$ for a decrease in $\tau$ and $\Omega_I=(1,1.1]$ for an increase.

\section{The effect of measurement error on the one-sided Synthetic-RZ control charts}
\label{sec:numerical}
This section reports how measurement error affects the two proposed charts. For convenience and without loss of generality, the in-control ratio and the target $ARL$ are fixed at $z_0=1$ and $ARL_0=200$, respectively, and the mean shift of $Y_{i,j}$ is set to $\delta_Y=1$. The remaining process parameters are allowed to vary over the following grid:
\begin{itemize}
\item $\gamma_X \in \{0.01, 0.2\}$, $\gamma_Y \in \{0.01, 0.2\}$,
\item $n \in \{1, 5, 7, 10, 15\}$,
\item $\rho_0 \in \{-0.8, -0.4, 0, 0.4, 0.8\}$,
\item $\tau \in \Omega_D = [0.9, 1)$ for the Synthetic-RZ$^-$-ME chart and $\tau \in \Omega_I = (1, 1.1]$ for the Synthetic-RZ$^+$-ME chart.
\end{itemize}

\subsection{Effect on the control limits}
Table \eqref{tab:zerostate:values} lists the optimal couples $(LCL^-, H^-)$ for the Synthetic-RZ$^-$-ME chart together with $(UCL^+, H^+)$ for the Synthetic-RZ$^+$-ME chart, computed in the zero state with measurement error present and with $\theta_X=\theta_Y=0.01$, $\eta_X=\eta_Y=0.28$, $m=1$, $\rho_M=0.5$ and $ARL_0=200$. The analogous steady-state quantities appear in Table \eqref{tab:steadystate:values}.

Two observations stand out. Although measurement error clearly alters the numerical values of the limits, the qualitative pattern matches what is seen when no error is present: both $LCL^-$ and $UCL^+$ vary with $n$ and with $\rho_0$, and the reciprocal relationship $LCL^- \approx 1/UCL^+$ continues to hold whenever $\gamma_X=\gamma_Y$. In addition, introducing measurement error pushes the two limits farther apart and modifies $H^-$ (respectively $H^+$); the larger the error, the wider the resulting interval between the limits.

\begin{table}[H]
\centering
\scalebox{0.78}{
\begin{tabular}{cccccccc}
\hline
$\gamma_X$ & $\gamma_Y$ & $\rho_0$ & $n=1$ & $n=5$ & $n=7$ & $n=10$ & $n=15$ \\
\hline
$0.01$ & $0.01$ & $-0.8$ & (0.9622, 14) & (0.9843, 6) & (0.9870, 5) & (0.9894, 4) & (0.9916, 3) \\
       &        &         & (1.0393, 14) & (1.0159, 6) & (1.0132, 5) & (1.0107, 4) & (1.0084, 3) \\
$0.01$ & $0.01$ & $-0.4$ & (0.9667, 13) & (0.9864, 5) & (0.9888, 4) & (0.9909, 3) & (0.9926, 3) \\
       &        &         & (1.0342, 12) & (1.0138, 5) & (1.0113, 4) & (1.0091, 3) & (1.0075, 3) \\
$0.01$ & $0.01$ & $0.0$ & (0.9721, 11) & (0.9887, 4) & (0.9908, 3) & (0.9923, 3) & (0.9940, 2) \\
       &        &         & (1.0287, 11) & (1.0114, 4) & (1.0093, 3) & (1.0078, 3) & (1.0060, 2) \\
$0.01$ & $0.01$ & $0.4$ & (0.9787, 8) & (0.9915, 3) & (0.9932, 2) & (0.9943, 2) & (0.9953, 2) \\
       &        &         & (1.0217, 8) & (1.0086, 3) & (1.0069, 2) & (1.0058, 2) & (1.0047, 2) \\
$0.01$ & $0.01$ & $0.8$ & (0.9876, 5) & (0.9950, 2) & (0.9958, 2) & (0.9965, 2) & (0.9974, 1) \\
       &        &         & (1.0126, 5) & (1.0050, 2) & (1.0042, 2) & (1.0035, 2) & (1.0026, 1) \\
\hline
$0.2$ & $0.2$ & $-0.8$ & (0.4704, 42) & (0.7280, 35) & (0.7662, 33) & (0.8015, 31) & (0.8354, 30) \\
       &        &         & (1.8806, 7) & (1.3540, 20) & (1.2924, 21) & (1.2399, 22) & (1.1916, 22) \\
$0.2$ & $0.2$ & $-0.4$ & (0.5129, 40) & (0.7559, 33) & (0.7905, 32) & (0.8226, 30) & (0.8536, 28) \\
       &        &         & (1.7672, 8) & (1.3092, 21) & (1.2549, 21) & (1.2096, 22) & (1.1678, 22) \\
$0.2$ & $0.2$ & $0.0$ & (0.5669, 38) & (0.7891, 31) & (0.8196, 30) & (0.8478, 28) & (0.8745, 27) \\
       &        &         & (1.6454, 10) & (1.2575, 21) & (1.2139, 22) & (1.1756, 22) & (1.1409, 22) \\
$0.2$ & $0.2$ & $0.4$ & (0.6428, 33) & (0.8312, 29) & (0.8562, 28) & (0.8792, 26) & (0.9011, 24) \\
       &        &         & (1.4918, 12) & (1.1971, 21) & (1.1643, 22) & (1.1353, 22) & (1.1085, 21) \\
$0.2$ & $0.2$ & $0.8$ & (0.7659, 27) & (0.8946, 24) & (0.9108, 23) & (0.9256, 21) & (0.9394, 19) \\
       &        &         & (1.2861, 15) & (1.1164, 21) & (1.0967, 20) & (1.0797, 19) & (1.0642, 18) \\
\hline
$0.01$ & $0.2$ & $-0.8$ & (0.7116, 42) & (0.8503, 32) & (0.8711, 30) & (0.8905, 28) & (0.9096, 25) \\
       &        &         & (1.4939, 6) & (1.1988, 18) & (1.1641, 19) & (1.1338, 19) & (1.1067, 19) \\
$0.01$ & $0.2$ & $-0.4$ & (0.7166, 42) & (0.8529, 32) & (0.8737, 29) & (0.8924, 28) & (0.9112, 25) \\
       &        &         & (1.4853, 6) & (1.1953, 18) & (1.1612, 19) & (1.1315, 19) & (1.1049, 19) \\
$0.01$ & $0.2$ & $0.0$ & (0.7217, 42) & (0.8558, 31) & (0.8760, 29) & (0.8946, 27) & (0.9128, 25) \\
       &        &         & (1.4765, 6) & (1.1918, 18) & (1.1583, 19) & (1.1291, 19) & (1.1030, 19) \\
$0.01$ & $0.2$ & $0.4$ & (0.7269, 42) & (0.8585, 31) & (0.8783, 29) & (0.8966, 27) & (0.9144, 25) \\
       &        &         & (1.4676, 6) & (1.1882, 18) & (1.1553, 19) & (1.1267, 19) & (1.1010, 19) \\
$0.01$ & $0.2$ & $0.8$ & (0.7323, 42) & (0.8613, 31) & (0.8807, 29) & (0.8986, 27) & (0.9163, 24) \\
       &        &         & (1.4585, 6) & (1.1846, 18) & (1.1523, 19) & (1.1242, 19) & (1.0991, 19) \\
\hline
$0.2$ & $0.01$ & $-0.8$ & (0.5455, 30) & (0.7967, 27) & (0.8285, 26) & (0.8564, 26) & (0.8834, 24) \\
       &        &         & (1.4683, 29) & (1.2057, 26) & (1.1736, 26) & (1.1446, 25) & (1.1176, 24) \\
$0.2$ & $0.01$ & $-0.4$ & (0.5511, 29) & (0.7995, 27) & (0.8310, 26) & (0.8590, 25) & (0.8852, 24) \\
       &        &         & (1.4565, 29) & (1.2020, 27) & (1.1701, 26) & (1.1418, 25) & (1.1153, 24) \\
$0.2$ & $0.01$ & $0.0$ & (0.5544, 30) & (0.8024, 27) & (0.8335, 26) & (0.8612, 25) & (0.8870, 24) \\
       &        &         & (1.4459, 30) & (1.1976, 27) & (1.1665, 26) & (1.1389, 25) & (1.1130, 24) \\
$0.2$ & $0.01$ & $0.4$ & (0.5603, 29) & (0.8054, 27) & (0.8361, 26) & (0.8634, 25) & (0.8893, 23) \\
       &        &         & (1.4338, 30) & (1.1932, 27) & (1.1629, 26) & (1.1359, 25) & (1.1106, 24) \\
$0.2$ & $0.01$ & $0.8$ & (0.5650, 29) & (0.8090, 26) & (0.8388, 26) & (0.8657, 25) & (0.8912, 23) \\
       &        &         & (1.4217, 30) & (1.1887, 27) & (1.1592, 26) & (1.1329, 25) & (1.1082, 24) \\
\hline
\end{tabular}}
\caption{Optimal couples $(LCL^-, H^-)$ (upper row) and $(UCL^+, H^+)$ (lower row) of the one-sided Synthetic-RZ-ME chart in the zero state, for $z_0=1$, $\theta_X=\theta_Y=0.01$, $\eta_X=\eta_Y=0.28$, $\rho_M=0.5$, $m=1$ and $ARL_0=200$.}
\label{tab:zerostate:values}
\end{table}

\input{tables/table_steadystate.tex}

\subsection{Effect on the chart performance}

To assess detection capability, Tables \eqref{tab:earl:zs:eq}--\eqref{tab:earl:ss:neq} collect the $EARL$ values obtained under two run-length scenarios---the zero-state case (Tables \eqref{tab:earl:zs:eq}--\eqref{tab:earl:zs:neq}) and the steady-state case (Tables \eqref{tab:earl:ss:eq}--\eqref{tab:earl:ss:neq})---and for two settings of the correlation, namely $\rho_0=\rho_1$ (Tables \eqref{tab:earl:zs:eq} and \eqref{tab:earl:ss:eq}) and $\rho_1=-0.8\neq\rho_0$ (Tables \eqref{tab:earl:zs:neq} and \eqref{tab:earl:ss:neq}).

The tabulated figures point to two findings of practical relevance. On the one hand, enlarging the sample from $n=1$ to $n=15$ lowers $EARL$ by a factor of roughly two to three for every parameter combination examined, so increasing the sample size remains the single most powerful lever for sharpening detection. On the other hand, the steady-state values exceed their zero-state counterparts by about $15$--$20\%$; this gap agrees with a familiar feature of synthetic schemes and provides a useful internal check on the Markov-chain computations.

\input{tables/table_earl_zs_eq.tex}
\input{tables/table_earl_zs_neq.tex}
\input{tables/table_earl_ss_eq.tex}
\input{tables/table_earl_ss_neq.tex}

\subsubsection{Effect of precision errors $\eta_X$ and $\eta_Y$}
\begin{figure}[H]
\centering
\includegraphics[width=0.95\textwidth]{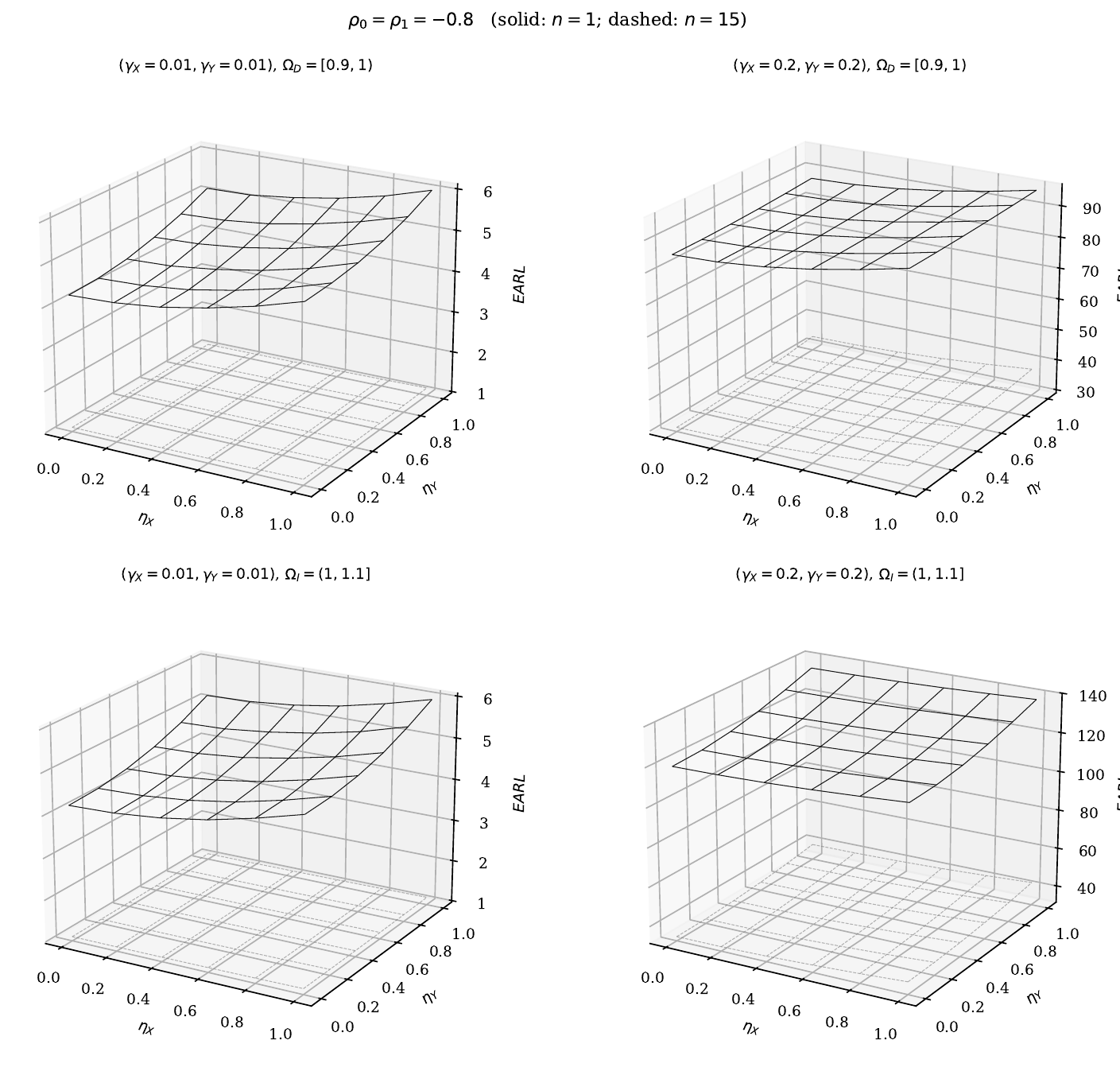}
\caption{Behaviour of $EARL$ as a function of the precision errors $\eta_X$ and $\eta_Y$, with $\theta_X=\theta_Y=0$, $\rho_M=0$, $n\in\{1,15\}$, $\gamma_X=\gamma_Y\in\{0.01,0.2\}$ and $\rho_0=\rho_1=-0.8$.}
\label{fig:1}
\end{figure}

\begin{figure}[H]
\centering
\includegraphics[width=0.95\textwidth]{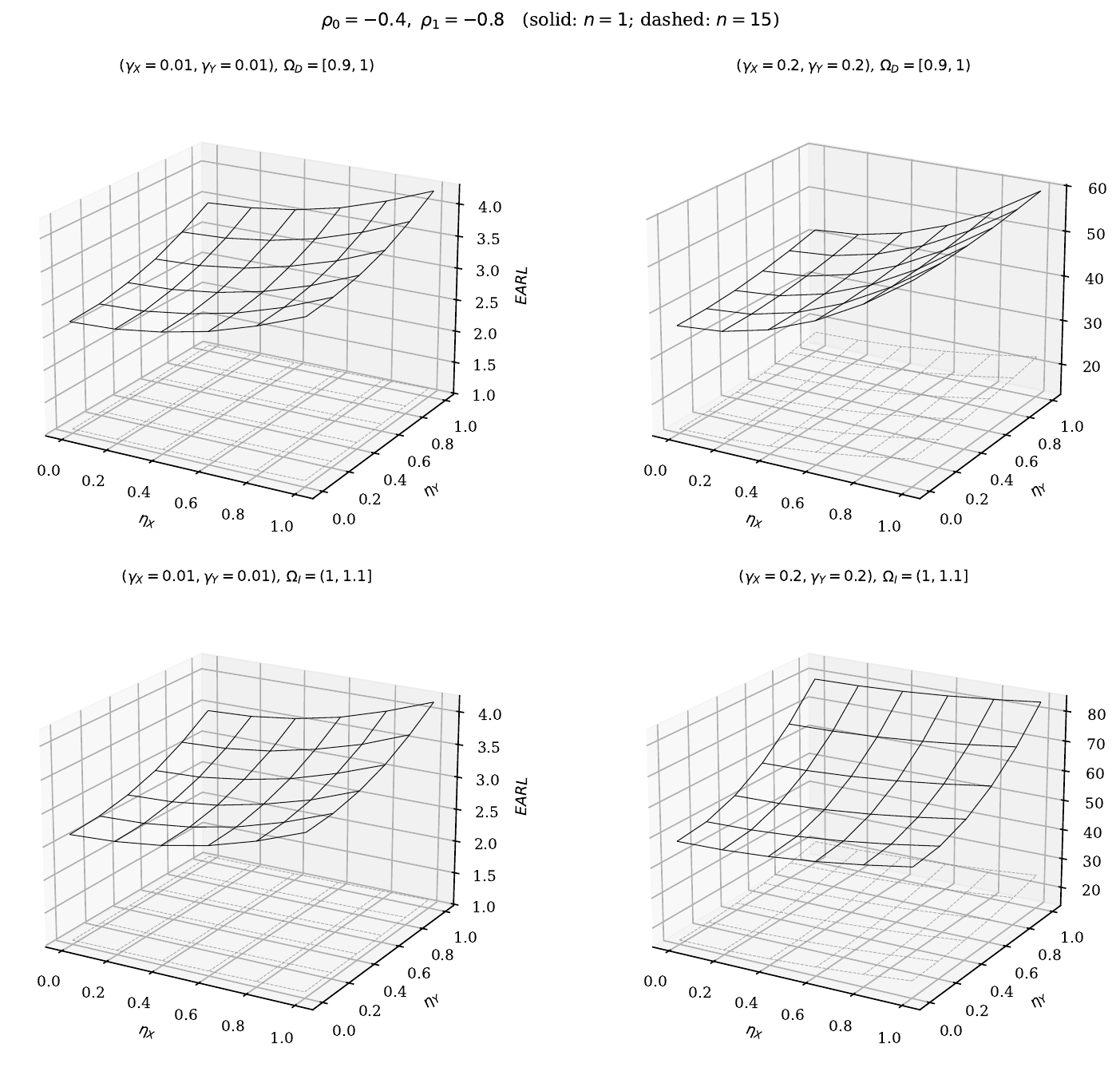}
\caption{Behaviour of $EARL$ versus $\eta_X$ and $\eta_Y$ under the same configuration, but with $\rho_0=-0.4$ and $\rho_1=-0.8$.}
\label{fig:2}
\end{figure}

We first turn to the precision component of the error, captured by $\eta_X$ and $\eta_Y$. Holding $\gamma_X = \gamma_Y \in \{0.01, 0.2\}$, $\rho_M=0$ and $\theta_X = \theta_Y = 0$, we treat two correlation configurations, $\rho_0 = \rho_1 = -0.8$ and $\rho_0 = -0.4,\,\rho_1 = -0.8$, and let each of $\eta_X$ and $\eta_Y$ range over $\{0, 0.1, 0.2, \ldots, 1\}$. The resulting $EARL$ curves are displayed in Figures \ref{fig:1}--\ref{fig:2}.

Because $EARL$ grows as $\eta_X$ and $\eta_Y$ increase, the precision error evidently harms the charts. To illustrate, with $n=1$, $\Omega_I=(1, 1.1]$, $\rho_0 = \rho_1 = -0.8$ and $\gamma_X = \gamma_Y = 0.2$, the chart yields $EARL=27.8$ at $\eta_X=\eta_Y=0.1$, and a still larger value once $\eta_X=\eta_Y=1$. The deterioration is nevertheless modest as long as the precision error stays small. Raising $n$ again helps substantially: across the whole range, the $EARL$ obtained for $n=15$ never exceeds that for $n=1$, whether or not measurement error is present.

\subsubsection{Effect of accuracy errors $\theta_X$ and $\theta_Y$}
\begin{figure}[H]
\centering
\includegraphics[width=0.95\textwidth]{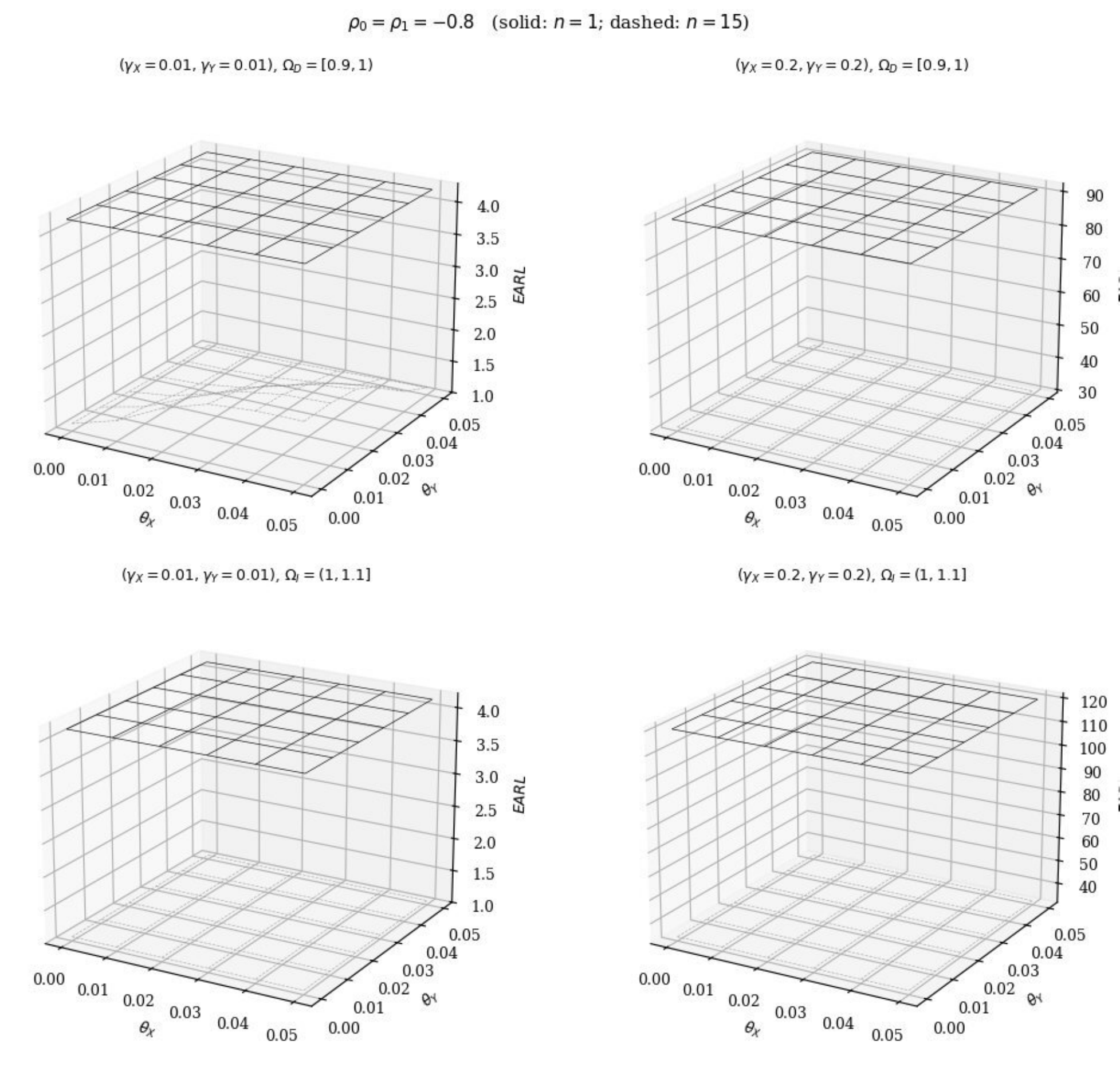}
\caption{Influence of the accuracy errors $\theta_X$ and $\theta_Y$ on $EARL$ for $\eta_X=\eta_Y=0$, $\rho_M=0$ and $\rho_0=\rho_1=-0.8$.}
\label{fig:3}
\end{figure}

\begin{figure}[H]
\centering
\includegraphics[width=0.95\textwidth]{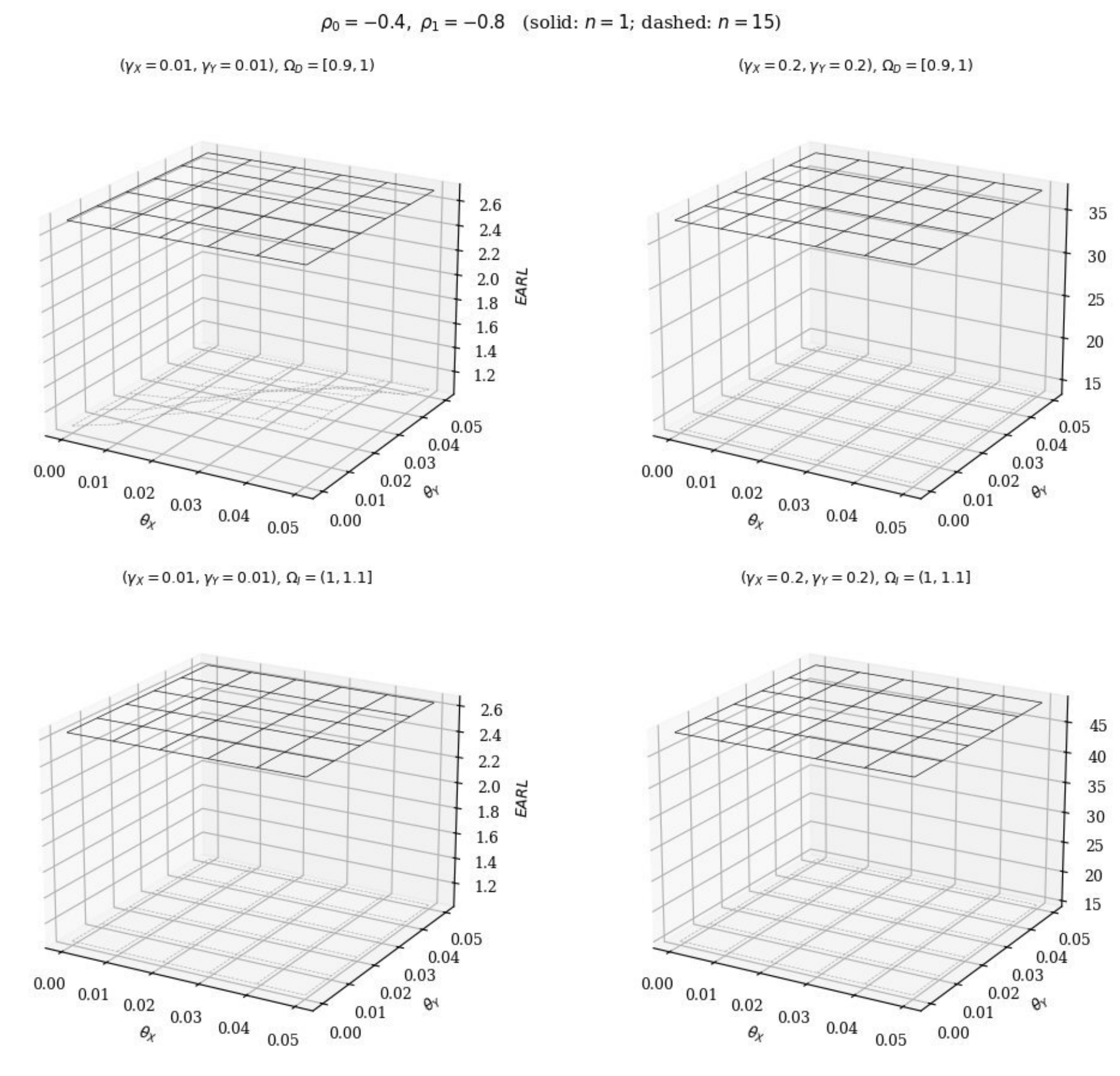}
\caption{Influence of $\theta_X$ and $\theta_Y$ on $EARL$ for $\rho_0=-0.4$ and $\rho_1=-0.8$.}
\label{fig:4}
\end{figure}

Next we consider the accuracy component, described by $\theta_X$ and $\theta_Y$. We let both quantities vary over $\{0, 0.01, 0.02, \ldots, 0.05\}$ while keeping $\gamma_X = \gamma_Y \in \{0.01, 0.2\}$, $\rho_M=0$ and $\eta_X = \eta_Y = 0$, again under $\rho_0 = \rho_1 = -0.8$ (Figure \ref{fig:3}) and $\rho_0 = -0.4,\,\rho_1 = -0.8$ (Figure \ref{fig:4}).

The figures reveal that how $\theta_X$ and $\theta_Y$ act on $EARL$ is governed by the values of $\rho_0$ and $\rho_1$. When $\rho_0 = \rho_1 = -0.8$, $EARL$ rises as $\theta_X$ increases and $\theta_Y$ decreases (Figure \ref{fig:3}). The opposite pattern emerges for $\rho_0 = -0.4,\,\rho_1 = -0.8$, where $EARL$ climbs when $\theta_X$ falls and $\theta_Y$ rises. As before, a larger $n$ markedly enhances performance, the $n=15$ curves lying everywhere below the $n=1$ curves irrespective of measurement error.

\subsubsection{Effect of the measurement-error correlation $\rho_M$}
\begin{figure}[H]
\centering
\includegraphics[width=0.95\textwidth]{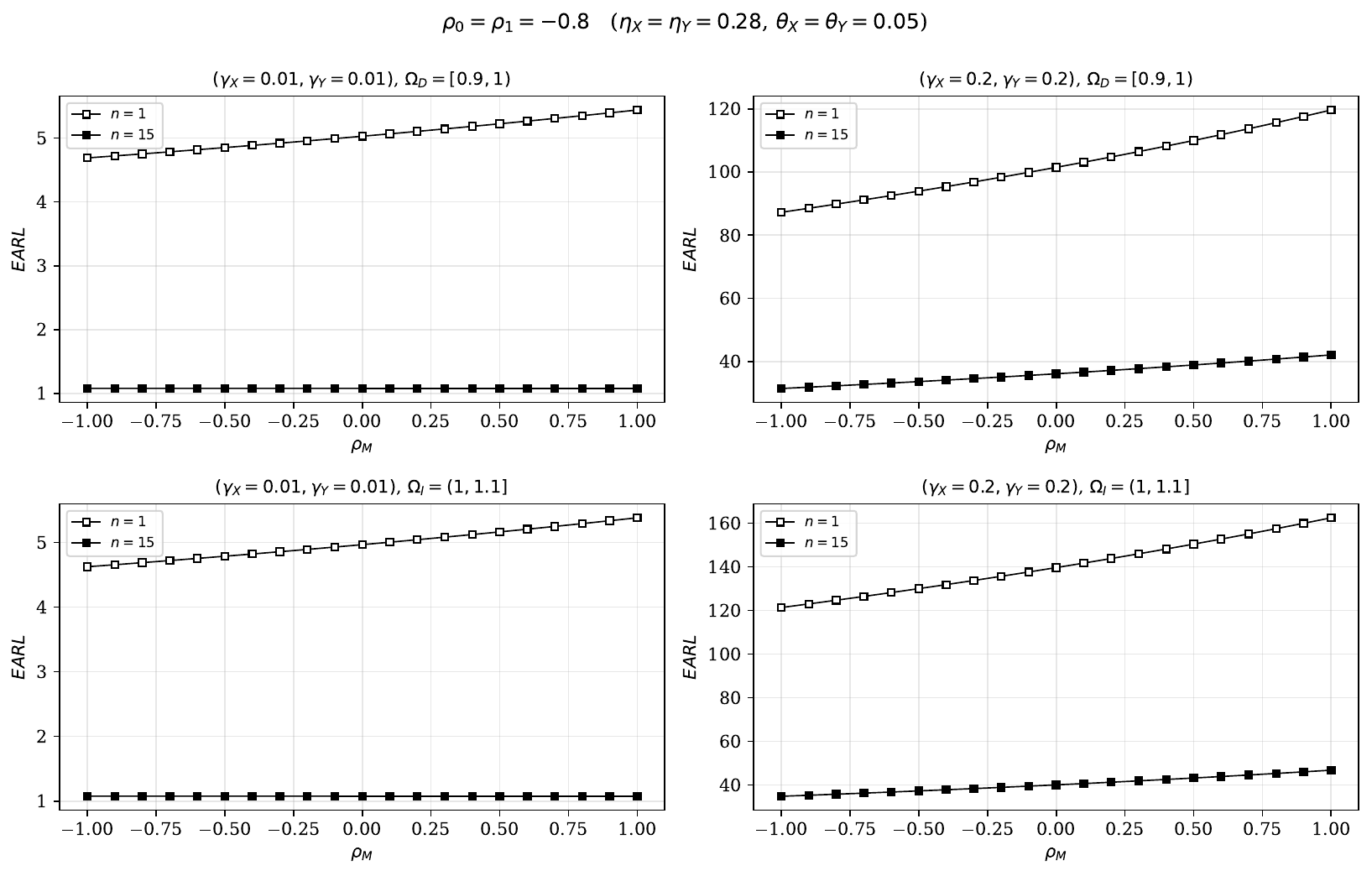}
\caption{Dependence of $EARL$ on $\rho_M$ for $\eta_X=\eta_Y=0.28$, $\theta_X=\theta_Y=0.05$ and $\rho_0=\rho_1=-0.8$.}
\label{fig:5}
\end{figure}

\begin{figure}[H]
\centering
\includegraphics[width=0.95\textwidth]{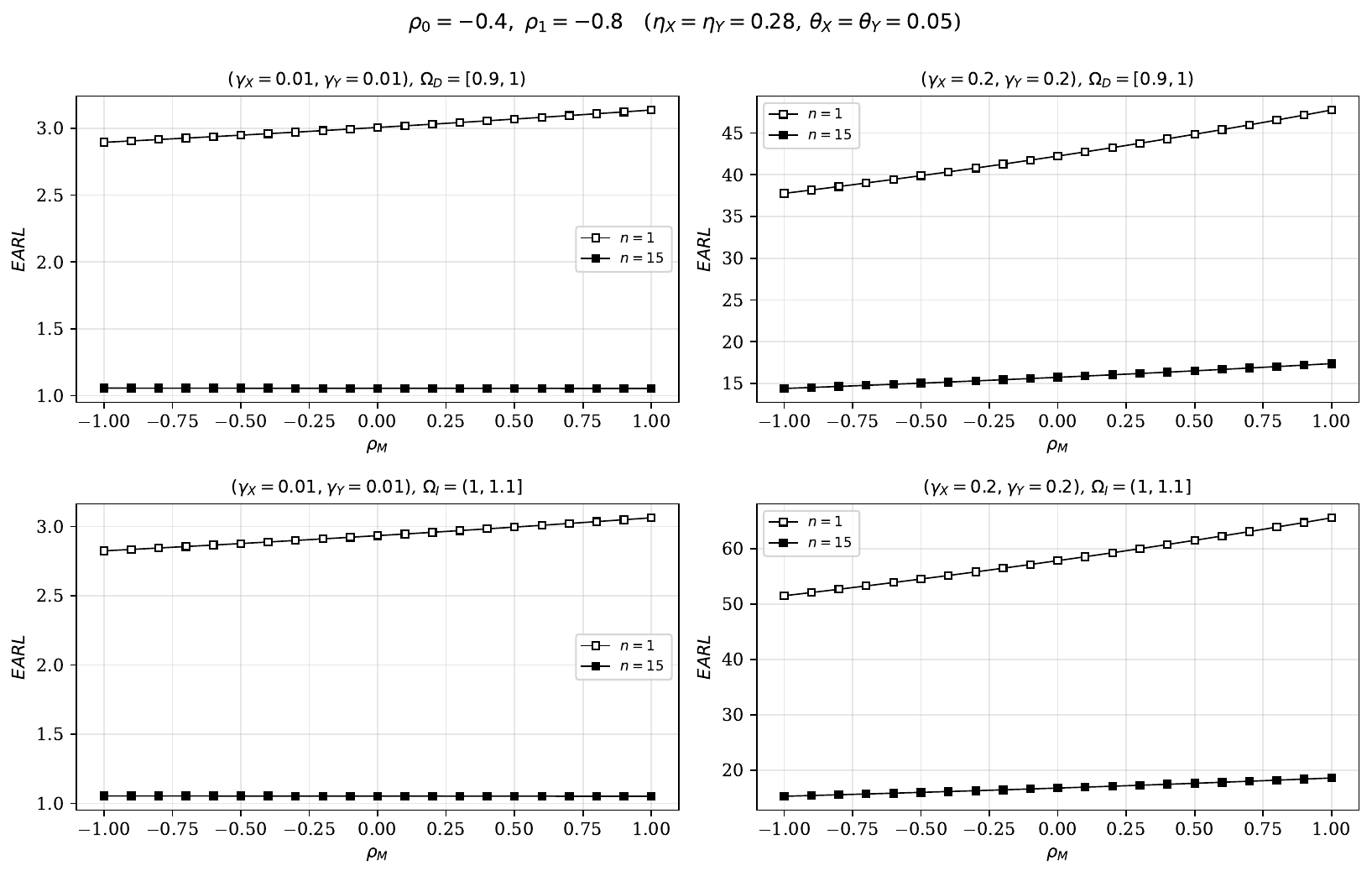}
\caption{Dependence of $EARL$ on $\rho_M$ for $\rho_0=-0.4$ and $\rho_1=-0.8$.}
\label{fig:6}
\end{figure}

The role of $\rho_M$ is depicted in Figure \ref{fig:5} ($\rho_0 = \rho_1 = -0.8$) and Figure \ref{fig:6} ($\rho_0 = -0.4,\,\rho_1 = -0.8$), obtained with $\gamma_X=\gamma_Y \in \{0.01, 0.2\}$, $\theta_X = \theta_Y = 0.05$, $\eta_X = \eta_Y = 0.28$ and $\rho_M$ ranging over $\{-1, -0.9, \ldots, 0.9, 1\}$. The curves indicate that $EARL$ responds smoothly to $\rho_M$. With both $\rho_0$ and $\rho_1$ negative, a higher $\rho_M$ produces a slight rise in $EARL$: a positive correlation between the measurement errors offsets part of the negative correlation in the process, so the standardized ratio becomes less concentrated and the chart loses a little of its sensitivity. The take-away is subtle yet practically relevant---the direction in which $\rho_M$ acts hinges on the signs of $\rho_0$ and $\rho_1$, and it would be a mistake to presume that mutually uncorrelated measurement errors invariably represent the least favourable situation.

\subsubsection{Effect of multiple measurements $m$}
\begin{figure}[H]
\centering
\includegraphics[width=0.95\textwidth]{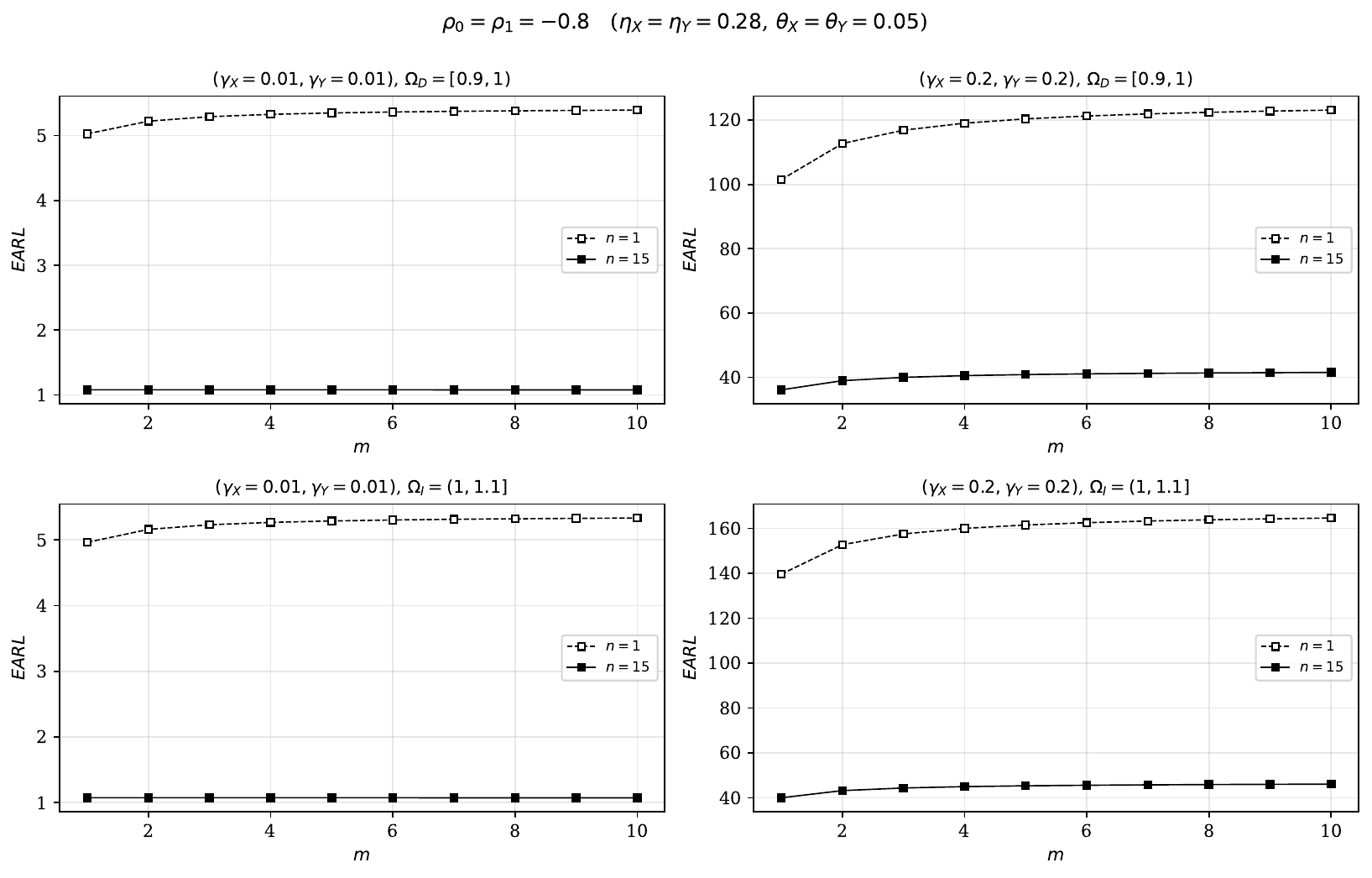}
\caption{Effect of the number of repeated measurements $m$ on $EARL$ for $\eta_X=\eta_Y=0.28$, $\theta_X=\theta_Y=0.05$ and $\rho_0=\rho_1=-0.8$.}
\label{fig:7}
\end{figure}

\begin{figure}[H]
\centering
\includegraphics[width=0.95\textwidth]{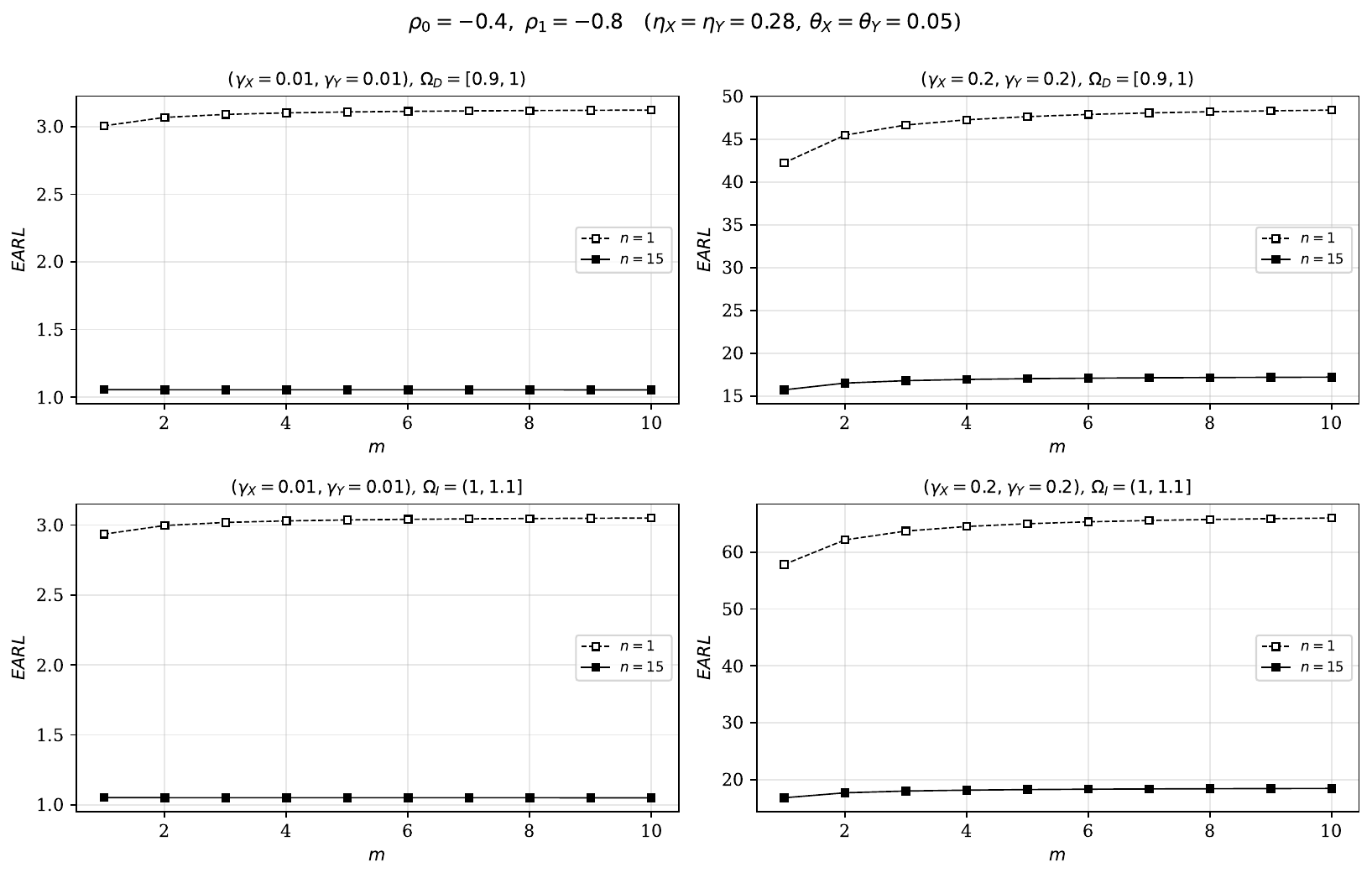}
\caption{Effect of $m$ on $EARL$ for $\rho_0=-0.4$ and $\rho_1=-0.8$.}
\label{fig:8}
\end{figure}

A recurring suggestion in the literature is that recording several measurements on each unit can offset the damage measurement error does to a chart. To test whether this holds for the Synthetic-RZ-ME chart, we evaluate $EARL$ as $m$ grows from $1$ to $10$, with $\gamma_X=\gamma_Y \in \{0.01, 0.2\}$, $\theta_X = \theta_Y = 0.05$ and $\eta_X = \eta_Y = 0.28$. Figures \ref{fig:7}--\ref{fig:8} report the outcome.

The conclusion is striking: enlarging $m$ brings no appreciable gain in the performance of the Synthetic-RZ-ME chart---the curves in the figures stay almost flat. This mirrors the result reported by \citet{Tran2016_Shewhart_RZ_ME} for the Shewhart-RZ chart. Consequently, repeated measurement of each item is not a worthwhile strategy for mitigating the adverse effect of measurement error on these charts. From a practical standpoint the message is valuable: effort is better spent on shrinking the intrinsic measurement error itself---that is, on lowering $\eta_X$, $\eta_Y$, $\theta_X$ and $\theta_Y$---than on accumulating repeated readings.

\subsection{Comparison with the Shewhart-RZ-ME chart}
To bring out the merit of the proposed charts, we set their performance against that of the Shewhart-RZ-ME chart studied by \citet{Tran2016_Shewhart_RZ_ME}. The comparison is admittedly not on an entirely level footing, since the two works rest on slightly different versions of the linear covariate error model and involve different chart families; even so, the Synthetic-RZ-ME charts are found to be markedly superior at picking up small-to-moderate shifts, particularly when the in-control configuration permits the Synthetic chart to operate with larger $H^-$ or $H^+$. This accords with the well-documented fact that Synthetic schemes usually detect shifts more effectively than Shewhart schemes. As an added benefit, the one-sided construction frees the Synthetic-RZ-ME charts from the $ARL$-biased behaviour that afflicts the two-sided Synthetic-RZ chart of \citet{Celano2016_Synthentic_RZ}.

\section{Illustrative example}
\label{sec:illustrative}
This section walks through an application of the one-sided Synthetic-RZ-ME chart under measurement error. We revisit the real waste-battery management problem at Italian battery-recycling plants previously considered by \citet{Tran2016_Shewhart_RZ_ME}.

At the outset of recycling, batteries arrive at collection sites and are placed in dedicated drums, sacks or boxes, collectively called ``batches''. Such batches often contain discarded items that are not recyclable batteries---small electronic devices, scrap metal and assorted waste---and removing this material is necessary to contain recycling costs. The plant monitors the ratio $Z$ of the recyclable-battery weight (denoted $X$) to the total batch weight (denoted $Y$) in order to quantify that cost. As in \citet{Tran2016_Shewhart_RZ_ME}, the in-control ratio of interest, beyond which an economic loss is incurred, is set at $z_0 = 0.95$.

Samples of size $n=5$ are taken at regular times $i=1,2,\ldots$, the batches having a nominal weight of $100$ kg. To reflect natural variability the batch weight is modelled as $Y \sim N(100, 1)$, and the recyclable-battery weight within a batch is likewise normal with mean $\mu_X = 95$ kg. At each time the sample mean weights $\bar{X}^*_i = \frac{1}{n}\sum_{j=1}^n X^*_{i,j}$ and $\bar{Y}^*_i = \frac{1}{n}\sum_{j=1}^n Y^*_{i,j}$ are recorded. The coefficients of variation are $\gamma_X=\gamma_Y=0.01$ and the in-control correlation is $\rho_0=0.8$.

\begin{figure}[H]
\centering
\includegraphics[width=0.95\textwidth]{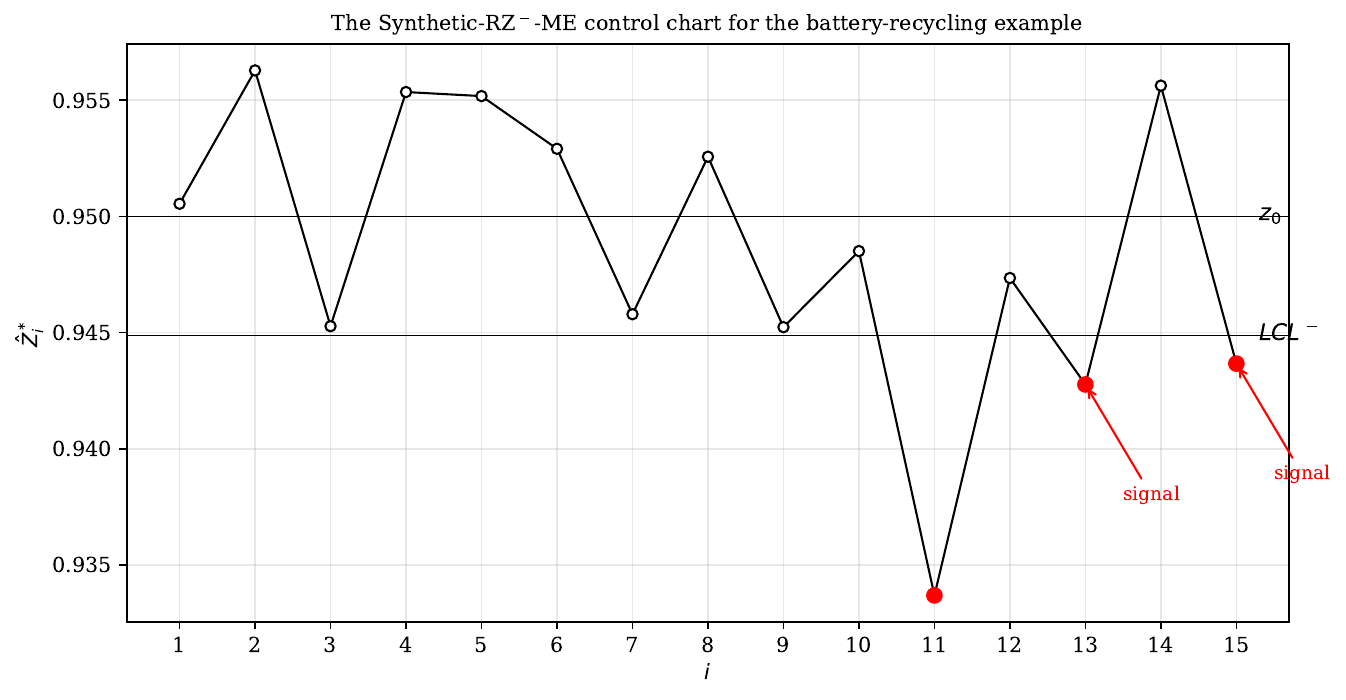}
\caption{Application of the Synthetic-RZ$^-$-ME chart, under measurement error, to the battery-recycling dataset.}
\label{fig:syn:example}
\end{figure}

The linear covariate error model is assigned the parameters $\theta_X = \theta_Y = 0$, $\eta_X = \eta_Y = 0.28$, $\rho_1 = 0.8$ and $\rho_M = 0$. From sample \#11 onward we simulate a downward shift amounting to at most $1\%$ of the in-control ratio $z_0$. For these values the optimal limit of the Synthetic-RZ$^-$-ME chart is obtained by the procedure of Section \ref{sec:implementation}. Figure \ref{fig:syn:example} shows the resulting chart. It signals an out-of-control condition at sample \#13---the plotted point falls below $LCL^-$ and the associated $CRL$ satisfies $CRL \leqslant H^-$---suggesting that an assignable cause has shifted the process out of control. The conclusion coincides with that of \citet{Tran2016_Shewhart_RZ_ME}.

\section{Concluding remarks}
\label{sec:conclusions}
This paper has examined how measurement error affects two one-sided Synthetic-RZ control charts used to monitor the ratio of two normal variables. The error was represented by a linear covariate error model, and we derived in full how the model parameters shift from the in-control to the out-of-control state without resorting to the restrictive assumption that the shift size is independent of the measurement error. The run-length properties of the Synthetic-RZ-ME charts were computed through a Markov chain approach covering both the zero-state and the steady-state $ARL$.

The numerical evidence confirms that measurement error degrades the performance of the one-sided Synthetic-RZ-ME charts. The precision errors ($\eta_X$ and $\eta_Y$) inflate the $EARL$, with larger errors producing worse detection. The accuracy errors ($\theta_X$ and $\theta_Y$) act in a more intricate fashion that is contingent on the values of $\rho_0$ and $\rho_1$. The correlation $\rho_M$ between the measurement errors of $X$ and $Y$ works in the opposite direction: a larger $\rho_M$ improves the chart. Enlarging the sample size $n$ enhances performance considerably in every case, irrespective of measurement error.

A noteworthy outcome is that collecting several measurements per item ($m>1$) does not effectively curb the harm measurement error does to the Synthetic-RZ-ME chart. This matches the conclusion reached for the Shewhart-RZ chart by \citet{Tran2016_Shewhart_RZ_ME} and constitutes a handy guideline for quality practitioners.

Several avenues remain open for further work. One is to attenuate the effect of measurement error by coupling the Synthetic-RZ chart with adaptive mechanisms such as variable sampling interval (VSI) or variable sample size (VSS), or with run rules. Another is to design other Synthetic-type charts---for instance Synthetic-EWMA-RZ or Synthetic-CUSUM-RZ---in the presence of measurement error. A final direction is to broaden the study to the monitoring of a multivariate ratio.

\bibliographystyle{plainnat}
\bibliography{paper}

\end{document}

%% file: tables/table_steadystate.tex
\begin{table}[H]
\centering
\scalebox{0.78}{
\begin{tabular}{cccccccc}
\hline
$\gamma_X$ & $\gamma_Y$ & $\rho_0$ & $n=1$ & $n=5$ & $n=7$ & $n=10$ & $n=15$ \\
\hline
$0.01$ & $0.01$ & $-0.8$ & (0.9665, 5) & (0.9858, 3) & (0.9880, 3) & (0.9899, 3) & (0.9922, 2) \\
       &        &         & (1.0346, 5) & (1.0144, 3) & (1.0122, 3) & (1.0102, 3) & (1.0079, 2) \\
$0.01$ & $0.01$ & $-0.4$ & (0.9703, 5) & (0.9874, 3) & (0.9893, 3) & (0.9915, 2) & (0.9931, 2) \\
       &        &         & (1.0306, 5) & (1.0128, 3) & (1.0108, 3) & (1.0085, 2) & (1.0070, 2) \\
$0.01$ & $0.01$ & $0.0$ & (0.9747, 5) & (0.9893, 3) & (0.9914, 2) & (0.9928, 2) & (0.9941, 2) \\
       &        &         & (1.0259, 5) & (1.0108, 3) & (1.0087, 2) & (1.0072, 2) & (1.0059, 2) \\
$0.01$ & $0.01$ & $0.4$ & (0.9806, 4) & (0.9920, 2) & (0.9933, 2) & (0.9944, 2) & (0.9950, 4) \\
       &        &         & (1.0198, 4) & (1.0080, 2) & (1.0068, 2) & (1.0057, 2) & (1.0051, 4) \\
$0.01$ & $0.01$ & $0.8$ & (0.9885, 3) & (0.9948, 3) & (0.9949, 13) & (0.9956, 17) & (0.9964, 18) \\
       &        &         & (1.0116, 3) & (1.0052, 3) & (1.0051, 12) & (1.0044, 17) & (1.0036, 18) \\
\hline
$0.2$ & $0.2$ & $-0.8$ & (0.5723, 3) & (0.7773, 4) & (0.8085, 4) & (0.8335, 5) & (0.8619, 5) \\
       &        &         & (1.6941, 2) & (1.2754, 3) & (1.2278, 3) & (1.1870, 3) & (1.1558, 4) \\
$0.2$ & $0.2$ & $-0.4$ & (0.6096, 3) & (0.8001, 4) & (0.8245, 5) & (0.8511, 5) & (0.8768, 5) \\
       &        &         & (1.5959, 2) & (1.2402, 3) & (1.1991, 3) & (1.1701, 4) & (1.1367, 4) \\
$0.2$ & $0.2$ & $0.0$ & (0.6559, 3) & (0.8273, 4) & (0.8487, 5) & (0.8720, 5) & (0.8943, 5) \\
       &        &         & (1.4892, 2) & (1.2009, 3) & (1.1734, 4) & (1.1429, 4) & (1.1151, 4) \\
$0.2$ & $0.2$ & $0.4$ & (0.7085, 4) & (0.8584, 5) & (0.8792, 5) & (0.8981, 5) & (0.9161, 5) \\
       &        &         & (1.3676, 2) & (1.1604, 4) & (1.1337, 4) & (1.1105, 4) & (1.0916, 5) \\
$0.2$ & $0.2$ & $0.8$ & (0.8094, 4) & (0.9108, 5) & (0.9242, 5) & (0.9351, 6) & (0.9468, 6) \\
       &        &         & (1.2257, 3) & (1.0954, 4) & (1.0820, 5) & (1.0680, 5) & (1.0551, 5) \\
\hline
$0.01$ & $0.2$ & $-0.8$ & (0.7543, 5) & (0.8736, 5) & (0.8911, 5) & (0.9056, 6) & (0.9216, 6) \\
       &        &         & (1.3976, 2) & (1.1560, 3) & (1.1340, 4) & (1.1098, 4) & (1.0879, 4) \\
$0.01$ & $0.2$ & $-0.4$ & (0.7586, 5) & (0.8758, 5) & (0.8911, 6) & (0.9072, 6) & (0.9230, 6) \\
       &        &         & (1.3906, 2) & (1.1533, 3) & (1.1317, 4) & (1.1079, 4) & (1.0864, 4) \\
$0.01$ & $0.2$ & $0.0$ & (0.7629, 5) & (0.8780, 5) & (0.8930, 6) & (0.9089, 6) & (0.9244, 6) \\
       &        &         & (1.3836, 2) & (1.1505, 3) & (1.1293, 4) & (1.1059, 4) & (1.0848, 4) \\
$0.01$ & $0.2$ & $0.4$ & (0.7674, 5) & (0.8803, 5) & (0.8950, 6) & (0.9106, 6) & (0.9258, 6) \\
       &        &         & (1.3764, 2) & (1.1477, 3) & (1.1268, 4) & (1.1039, 4) & (1.0832, 4) \\
$0.01$ & $0.2$ & $0.8$ & (0.7719, 5) & (0.8826, 5) & (0.8971, 6) & (0.9124, 6) & (0.9272, 6) \\
       &        &         & (1.3691, 2) & (1.1448, 3) & (1.1244, 4) & (1.1019, 4) & (1.0816, 4) \\
\hline
$0.2$ & $0.01$ & $-0.8$ & (0.6448, 3) & (0.8344, 4) & (0.8599, 4) & (0.8797, 5) & (0.9017, 5) \\
       &        &         & (1.3643, 3) & (1.1675, 4) & (1.1415, 4) & (1.1213, 5) & (1.0990, 5) \\
$0.2$ & $0.01$ & $-0.4$ & (0.6489, 3) & (0.8369, 4) & (0.8620, 4) & (0.8816, 5) & (0.9033, 5) \\
       &        &         & (1.3558, 3) & (1.1641, 4) & (1.1387, 4) & (1.1190, 5) & (1.0971, 5) \\
$0.2$ & $0.01$ & $0.0$ & (0.6530, 3) & (0.8393, 4) & (0.8642, 4) & (0.8835, 5) & (0.9048, 5) \\
       &        &         & (1.3471, 3) & (1.1607, 4) & (1.1358, 4) & (1.1166, 5) & (1.0952, 5) \\
$0.2$ & $0.01$ & $0.4$ & (0.6573, 3) & (0.8419, 4) & (0.8664, 4) & (0.8854, 5) & (0.9065, 5) \\
       &        &         & (1.3384, 3) & (1.1572, 4) & (1.1329, 4) & (1.1141, 5) & (1.0932, 5) \\
$0.2$ & $0.01$ & $0.8$ & (0.6617, 3) & (0.8444, 4) & (0.8687, 4) & (0.8874, 5) & (0.9081, 5) \\
       &        &         & (1.3296, 3) & (1.1537, 4) & (1.1333, 5) & (1.1116, 5) & (1.0912, 5) \\
\hline
\end{tabular}}
\caption{Values of $(LCL^-, H^-)$ (first row) and $(UCL^+, H^+)$ (second row) for the one-sided Synthetic-RZ-ME chart in steady state, with $z_0=1$, $\theta_X=\theta_Y=0.01$, $\eta_X=\eta_Y=0.28$, $\rho_M=0.5$, $m=1$, $ARL_0=200$.}
\label{tab:steadystate:values}
\end{table}

%% file: tables/table_earl_zs_eq.tex
\begin{table}[H]
\centering
\scalebox{0.85}{
\begin{tabular}{cccccccc}
\hline
$\gamma_X$ & $\gamma_Y$ & $\rho_0$ & $n=1$ & $n=5$ & $n=7$ & $n=10$ & $n=15$ \\
\hline
$0.01$ & $0.01$ & $-0.8$ & 4.3 & 1.5 & 1.3 & 1.1 & 1.1 \\
       &        &         & 4.3 & 1.4 & 1.3 & 1.1 & 1.1 \\
$0.01$ & $0.01$ & $-0.4$ & 3.6 & 1.3 & 1.2 & 1.1 & 1.0 \\
       &        &         & 3.5 & 1.3 & 1.2 & 1.1 & 1.0 \\
$0.01$ & $0.01$ & $0.0$ & 2.8 & 1.2 & 1.1 & 1.1 & 1.0 \\
       &        &         & 2.7 & 1.2 & 1.1 & 1.0 & 1.0 \\
$0.01$ & $0.01$ & $0.4$ & 2.0 & 1.1 & 1.0 & 1.0 & 1.0 \\
       &        &         & 1.9 & 1.1 & 1.0 & 1.0 & 1.0 \\
$0.01$ & $0.01$ & $0.8$ & 1.2 & 1.0 & 1.0 & 1.0 & 1.0 \\
       &        &         & 1.2 & 1.0 & 1.0 & 1.0 & 1.0 \\
\hline
$0.2$ & $0.2$ & $-0.8$ & 85.9 & 50.8 & 44.0 & 37.4 & 30.7 \\
       &        &         & 120.5 & 61.4 & 51.6 & 42.6 & 33.8 \\
$0.2$ & $0.2$ & $-0.4$ & 81.2 & 45.9 & 39.5 & 33.3 & 27.0 \\
       &        &         & 111.7 & 54.4 & 45.4 & 37.2 & 29.3 \\
$0.2$ & $0.2$ & $0.0$ & 74.7 & 39.8 & 33.9 & 28.2 & 22.6 \\
       &        &         & 99.9 & 45.8 & 38.0 & 30.8 & 24.0 \\
$0.2$ & $0.2$ & $0.4$ & 64.5 & 31.6 & 26.4 & 21.6 & 16.9 \\
       &        &         & 82.6 & 35.0 & 28.6 & 22.9 & 17.5 \\
$0.2$ & $0.2$ & $0.8$ & 44.4 & 18.3 & 14.8 & 11.6 & 8.7 \\
       &        &         & 52.0 & 19.2 & 15.2 & 11.8 & 8.7 \\
\hline
$0.01$ & $0.2$ & $-0.8$ & 53.4 & 27.6 & 23.2 & 19.1 & 15.0 \\
       &        &         & 88.6 & 35.7 & 28.9 & 22.8 & 17.3 \\
$0.01$ & $0.2$ & $-0.4$ & 52.5 & 27.0 & 22.7 & 18.6 & 14.6 \\
       &        &         & 87.6 & 35.0 & 28.3 & 22.4 & 16.9 \\
$0.01$ & $0.2$ & $0.0$ & 51.7 & 26.4 & 22.2 & 18.2 & 14.3 \\
       &        &         & 86.5 & 34.3 & 27.7 & 21.9 & 16.5 \\
$0.01$ & $0.2$ & $0.4$ & 50.8 & 25.9 & 21.7 & 17.8 & 13.9 \\
       &        &         & 85.4 & 33.7 & 27.1 & 21.4 & 16.1 \\
$0.01$ & $0.2$ & $0.8$ & 49.9 & 25.3 & 21.2 & 17.3 & 13.5 \\
       &        &         & 84.2 & 33.0 & 26.5 & 20.9 & 15.7 \\
\hline
$0.2$ & $0.01$ & $-0.8$ & 78.1 & 38.6 & 32.3 & 26.5 & 20.8 \\
       &        &         & 75.3 & 36.2 & 30.2 & 24.7 & 19.3 \\
$0.2$ & $0.01$ & $-0.4$ & 77.5 & 38.1 & 31.8 & 26.0 & 20.4 \\
       &        &         & 73.8 & 35.5 & 29.5 & 24.1 & 18.8 \\
$0.2$ & $0.01$ & $0.0$ & 77.0 & 37.5 & 31.3 & 25.5 & 20.0 \\
       &        &         & 72.4 & 34.7 & 28.9 & 23.5 & 18.4 \\
$0.2$ & $0.01$ & $0.4$ & 76.4 & 36.9 & 30.8 & 25.1 & 19.5 \\
       &        &         & 70.9 & 33.8 & 28.2 & 22.9 & 17.9 \\
$0.2$ & $0.01$ & $0.8$ & 75.7 & 36.3 & 30.2 & 24.6 & 19.1 \\
       &        &         & 69.4 & 33.0 & 27.4 & 22.3 & 17.4 \\
\hline
\end{tabular}}
\caption{$EARL$ of the Synthetic-RZ$^-$-ME chart with $\tau \in \Omega_D=[0.9,1)$ (first row) and Synthetic-RZ$^+$-ME chart with $\tau \in \Omega_I=(1,1.1]$ (second row) for $\rho_0=\rho_1$, zero state.}
\label{tab:earl:zs:eq}
\end{table}

%% file: tables/table_earl_zs_neq.tex
\begin{table}[H]
\centering
\scalebox{0.85}{
\begin{tabular}{cccccccc}
\hline
$\gamma_X$ & $\gamma_Y$ & $\rho_0$ & $n=1$ & $n=5$ & $n=7$ & $n=10$ & $n=15$ \\
\hline
$0.01$ & $0.01$ & $-0.8$ & 4.3 & 1.5 & 1.3 & 1.1 & 1.1 \\
       &        &         & 4.3 & 1.4 & 1.3 & 1.1 & 1.1 \\
$0.01$ & $0.01$ & $-0.4$ & 2.8 & 1.3 & 1.2 & 1.1 & 1.0 \\
       &        &         & 2.7 & 1.3 & 1.2 & 1.1 & 1.0 \\
$0.01$ & $0.01$ & $0.0$ & 1.9 & 1.2 & 1.1 & 1.1 & 1.0 \\
       &        &         & 1.9 & 1.2 & 1.1 & 1.1 & 1.0 \\
$0.01$ & $0.01$ & $0.4$ & 1.5 & 1.1 & 1.1 & 1.0 & 1.0 \\
       &        &         & 1.5 & 1.1 & 1.1 & 1.0 & 1.0 \\
$0.01$ & $0.01$ & $0.8$ & 1.2 & 1.1 & 1.0 & 1.0 & 1.0 \\
       &        &         & 1.2 & 1.1 & 1.0 & 1.0 & 1.0 \\
\hline
$0.2$ & $0.2$ & $-0.8$ & 85.9 & 50.8 & 44.0 & 37.4 & 30.7 \\
       &        &         & 120.5 & 61.4 & 51.6 & 42.6 & 33.8 \\
$0.2$ & $0.2$ & $-0.4$ & 38.1 & 22.4 & 19.7 & 17.1 & 14.3 \\
       &        &         & 52.9 & 25.7 & 22.0 & 18.6 & 15.2 \\
$0.2$ & $0.2$ & $0.0$ & 17.6 & 10.9 & 9.8 & 8.6 & 7.4 \\
       &        &         & 21.8 & 11.6 & 10.3 & 9.0 & 7.6 \\
$0.2$ & $0.2$ & $0.4$ & 8.4 & 5.7 & 5.2 & 4.7 & 4.1 \\
       &        &         & 8.9 & 5.8 & 5.3 & 4.7 & 4.1 \\
$0.2$ & $0.2$ & $0.8$ & 3.9 & 3.0 & 2.8 & 2.6 & 2.3 \\
       &        &         & 3.8 & 3.0 & 2.8 & 2.5 & 2.3 \\
\hline
$0.01$ & $0.2$ & $-0.8$ & 53.4 & 27.6 & 23.2 & 19.1 & 15.0 \\
       &        &         & 88.6 & 35.7 & 28.9 & 22.8 & 17.3 \\
$0.01$ & $0.2$ & $-0.4$ & 44.2 & 23.8 & 20.2 & 16.7 & 13.2 \\
       &        &         & 80.7 & 31.8 & 25.7 & 20.4 & 15.5 \\
$0.01$ & $0.2$ & $0.0$ & 36.9 & 20.7 & 17.6 & 14.7 & 11.7 \\
       &        &         & 73.4 & 28.4 & 23.0 & 18.2 & 13.9 \\
$0.01$ & $0.2$ & $0.4$ & 31.1 & 18.0 & 15.5 & 13.0 & 10.4 \\
       &        &         & 66.8 & 25.3 & 20.5 & 16.3 & 12.5 \\
$0.01$ & $0.2$ & $0.8$ & 26.4 & 15.8 & 13.6 & 11.5 & 9.3 \\
       &        &         & 60.6 & 22.6 & 18.4 & 14.6 & 11.2 \\
\hline
$0.2$ & $0.01$ & $-0.8$ & 78.1 & 38.6 & 32.3 & 26.5 & 20.8 \\
       &        &         & 75.3 & 36.2 & 30.2 & 24.7 & 19.3 \\
$0.2$ & $0.01$ & $-0.4$ & 71.6 & 34.6 & 28.9 & 23.7 & 18.6 \\
       &        &         & 60.9 & 30.9 & 26.0 & 21.4 & 16.9 \\
$0.2$ & $0.01$ & $0.0$ & 65.7 & 31.0 & 25.9 & 21.2 & 16.7 \\
       &        &         & 49.6 & 26.5 & 22.4 & 18.6 & 14.8 \\
$0.2$ & $0.01$ & $0.4$ & 60.2 & 27.9 & 23.3 & 19.1 & 15.0 \\
       &        &         & 40.7 & 22.8 & 19.4 & 16.2 & 13.0 \\
$0.2$ & $0.01$ & $0.8$ & 55.1 & 25.0 & 20.9 & 17.1 & 13.5 \\
       &        &         & 33.7 & 19.7 & 16.9 & 14.2 & 11.5 \\
\hline
\end{tabular}}
\caption{$EARL$ for $\rho_1=-0.8 \neq \rho_0$, zero state.}
\label{tab:earl:zs:neq}
\end{table}

%% file: tables/table_earl_ss_eq.tex
\begin{table}[H]
\centering
\scalebox{0.85}{
\begin{tabular}{cccccccc}
\hline
$\gamma_X$ & $\gamma_Y$ & $\rho_0$ & $n=1$ & $n=5$ & $n=7$ & $n=10$ & $n=15$ \\
\hline
$0.01$ & $0.01$ & $-0.8$ & 6.2 & 2.6 & 2.3 & 2.1 & 2.0 \\
       &        &         & 6.1 & 2.5 & 2.3 & 2.1 & 2.0 \\
$0.01$ & $0.01$ & $-0.4$ & 5.3 & 2.4 & 2.2 & 2.1 & 2.0 \\
       &        &         & 5.2 & 2.3 & 2.1 & 2.0 & 2.0 \\
$0.01$ & $0.01$ & $0.0$ & 4.3 & 2.2 & 2.1 & 2.0 & 1.9 \\
       &        &         & 4.2 & 2.2 & 2.1 & 2.0 & 1.9 \\
$0.01$ & $0.01$ & $0.4$ & 3.3 & 2.0 & 2.0 & 1.9 & 1.9 \\
       &        &         & 3.2 & 2.0 & 2.0 & 1.9 & 1.9 \\
$0.01$ & $0.01$ & $0.8$ & 2.2 & 1.9 & 1.7 & 1.7 & 1.7 \\
       &        &         & 2.2 & 1.9 & 1.8 & 1.7 & 1.7 \\
\hline
$0.2$ & $0.2$ & $-0.8$ & 98.1 & 59.0 & 51.6 & 43.9 & 36.3 \\
       &        &         & 124.5 & 68.2 & 58.3 & 48.7 & 39.1 \\
$0.2$ & $0.2$ & $-0.4$ & 92.8 & 53.6 & 46.2 & 39.2 & 32.1 \\
       &        &         & 116.3 & 61.0 & 51.6 & 42.7 & 34.1 \\
$0.2$ & $0.2$ & $0.0$ & 85.5 & 46.8 & 39.9 & 33.5 & 27.1 \\
       &        &         & 105.2 & 52.1 & 43.5 & 35.8 & 28.3 \\
$0.2$ & $0.2$ & $0.4$ & 73.5 & 37.3 & 31.4 & 26.0 & 20.7 \\
       &        &         & 88.6 & 40.3 & 33.3 & 27.1 & 21.1 \\
$0.2$ & $0.2$ & $0.8$ & 51.3 & 22.2 & 18.2 & 14.6 & 11.3 \\
       &        &         & 57.8 & 22.9 & 18.5 & 14.7 & 11.2 \\
\hline
$0.01$ & $0.2$ & $-0.8$ & 62.6 & 33.0 & 28.0 & 23.2 & 18.6 \\
       &        &         & 93.2 & 40.8 & 33.4 & 26.9 & 20.8 \\
$0.01$ & $0.2$ & $-0.4$ & 61.7 & 32.4 & 27.4 & 22.8 & 18.2 \\
       &        &         & 92.2 & 40.1 & 32.8 & 26.3 & 20.3 \\
$0.01$ & $0.2$ & $0.0$ & 60.8 & 31.8 & 26.9 & 22.3 & 17.8 \\
       &        &         & 91.1 & 39.4 & 32.2 & 25.8 & 19.9 \\
$0.01$ & $0.2$ & $0.4$ & 59.8 & 31.1 & 26.3 & 21.8 & 17.4 \\
       &        &         & 90.0 & 38.6 & 31.5 & 25.3 & 19.5 \\
$0.01$ & $0.2$ & $0.8$ & 58.8 & 30.5 & 25.7 & 21.3 & 16.9 \\
       &        &         & 88.9 & 37.9 & 30.9 & 24.7 & 19.0 \\
\hline
$0.2$ & $0.01$ & $-0.8$ & 87.4 & 45.0 & 38.0 & 31.4 & 25.0 \\
       &        &         & 84.1 & 42.1 & 35.5 & 29.2 & 23.2 \\
$0.2$ & $0.01$ & $-0.4$ & 86.8 & 44.3 & 37.4 & 30.8 & 24.5 \\
       &        &         & 82.7 & 41.3 & 34.7 & 28.6 & 22.7 \\
$0.2$ & $0.01$ & $0.0$ & 86.1 & 43.7 & 36.8 & 30.3 & 24.0 \\
       &        &         & 81.3 & 40.4 & 34.0 & 27.9 & 22.1 \\
$0.2$ & $0.01$ & $0.4$ & 85.3 & 43.0 & 36.2 & 29.7 & 23.6 \\
       &        &         & 79.8 & 39.5 & 33.2 & 27.3 & 21.6 \\
$0.2$ & $0.01$ & $0.8$ & 84.6 & 42.3 & 35.5 & 29.2 & 23.1 \\
       &        &         & 78.3 & 38.6 & 32.3 & 26.6 & 21.1 \\
\hline
\end{tabular}}
\caption{$EARL$ for $\rho_0=\rho_1$, steady state.}
\label{tab:earl:ss:eq}
\end{table}

%% file: tables/table_earl_ss_neq.tex
\begin{table}[H]
\centering
\scalebox{0.85}{
\begin{tabular}{cccccccc}
\hline
$\gamma_X$ & $\gamma_Y$ & $\rho_0$ & $n=1$ & $n=5$ & $n=7$ & $n=10$ & $n=15$ \\
\hline
$0.01$ & $0.01$ & $-0.8$ & 6.2 & 2.6 & 2.3 & 2.1 & 2.0 \\
       &        &         & 6.1 & 2.5 & 2.3 & 2.1 & 2.0 \\
$0.01$ & $0.01$ & $-0.4$ & 4.4 & 2.3 & 2.2 & 2.1 & 2.0 \\
       &        &         & 4.3 & 2.3 & 2.1 & 2.1 & 2.0 \\
$0.01$ & $0.01$ & $0.0$ & 3.3 & 2.2 & 2.1 & 2.0 & 2.0 \\
       &        &         & 3.3 & 2.2 & 2.1 & 2.0 & 1.9 \\
$0.01$ & $0.01$ & $0.4$ & 2.7 & 2.1 & 2.0 & 2.0 & 1.9 \\
       &        &         & 2.7 & 2.1 & 2.0 & 2.0 & 1.9 \\
$0.01$ & $0.01$ & $0.8$ & 2.3 & 2.0 & 1.8 & 1.7 & 1.7 \\
       &        &         & 2.3 & 2.0 & 1.8 & 1.7 & 1.7 \\
\hline
$0.2$ & $0.2$ & $-0.8$ & 98.1 & 59.0 & 51.6 & 43.9 & 36.3 \\
       &        &         & 124.5 & 68.2 & 58.3 & 48.7 & 39.1 \\
$0.2$ & $0.2$ & $-0.4$ & 53.9 & 32.0 & 27.7 & 24.0 & 20.2 \\
       &        &         & 67.5 & 36.5 & 31.6 & 26.2 & 21.5 \\
$0.2$ & $0.2$ & $0.0$ & 28.9 & 17.6 & 15.3 & 13.5 & 11.6 \\
       &        &         & 35.5 & 19.7 & 16.6 & 14.4 & 12.2 \\
$0.2$ & $0.2$ & $0.4$ & 14.3 & 9.5 & 8.7 & 7.8 & 6.9 \\
       &        &         & 18.0 & 10.1 & 9.1 & 8.1 & 6.9 \\
$0.2$ & $0.2$ & $0.8$ & 6.9 & 5.2 & 4.8 & 4.4 & 4.1 \\
       &        &         & 7.5 & 5.3 & 4.8 & 4.5 & 4.1 \\
\hline
$0.01$ & $0.2$ & $-0.8$ & 62.6 & 33.0 & 28.0 & 23.2 & 18.6 \\
       &        &         & 93.2 & 40.8 & 33.4 & 26.9 & 20.8 \\
$0.01$ & $0.2$ & $-0.4$ & 54.9 & 29.7 & 25.2 & 21.1 & 16.9 \\
       &        &         & 86.5 & 37.6 & 30.7 & 24.7 & 19.1 \\
$0.01$ & $0.2$ & $0.0$ & 48.3 & 26.8 & 22.8 & 19.1 & 15.5 \\
       &        &         & 80.3 & 34.5 & 28.2 & 22.7 & 17.6 \\
$0.01$ & $0.2$ & $0.4$ & 42.5 & 24.2 & 20.6 & 17.4 & 14.1 \\
       &        &         & 74.5 & 31.8 & 25.9 & 20.9 & 16.3 \\
$0.01$ & $0.2$ & $0.8$ & 37.6 & 21.9 & 18.7 & 15.9 & 13.0 \\
       &        &         & 69.0 & 29.2 & 23.8 & 19.2 & 15.0 \\
\hline
$0.2$ & $0.01$ & $-0.8$ & 87.4 & 45.0 & 38.0 & 31.4 & 25.0 \\
       &        &         & 84.1 & 42.1 & 35.5 & 29.2 & 23.2 \\
$0.2$ & $0.01$ & $-0.4$ & 81.6 & 41.4 & 35.0 & 28.8 & 23.0 \\
       &        &         & 73.4 & 37.7 & 31.9 & 26.3 & 21.0 \\
$0.2$ & $0.01$ & $0.0$ & 76.2 & 38.1 & 32.2 & 26.5 & 21.2 \\
       &        &         & 64.2 & 33.7 & 28.7 & 23.7 & 19.1 \\
$0.2$ & $0.01$ & $0.4$ & 71.0 & 35.1 & 29.6 & 24.4 & 19.5 \\
       &        &         & 56.3 & 30.3 & 25.9 & 21.4 & 17.3 \\
$0.2$ & $0.01$ & $0.8$ & 66.2 & 32.3 & 27.3 & 22.4 & 18.0 \\
       &        &         & 49.4 & 27.2 & 23.0 & 19.4 & 15.8 \\
\hline
\end{tabular}}
\caption{$EARL$ for $\rho_1=-0.8 \neq \rho_0$, steady state.}
\label{tab:earl:ss:neq}
\end{table}